\documentclass[12pt]{article}
\usepackage{amsmath,amssymb,amsthm}
\usepackage[margin=1in]{geometry}
\usepackage{hyperref}
\usepackage{tikz}
\usepackage{graphicx}
\usepackage{url}
\usepackage[mathlines]{lineno}
\usepackage{multirow}
\usepackage{dsfont} % needed for double-struck 1
\usepackage{color}
\definecolor{red}{rgb}{1,0,0}

\definecolor{blue}{rgb}{0,0,1}

\definecolor{green}{rgb}{0,.6,0}

\usepackage{tikz}
\usetikzlibrary{positioning}
\tikzset{
every picture/.style=thick,
bluenode/.style={circle, draw=black, fill=blue!40, very thick, minimum size=6mm},
whitenode/.style={circle, draw=black, fill=black!5, very thick, minimum size=2mm},
squarednode/.style={rectangle, draw=red!60, fill=red!5, very thick, minimum size=5mm},
every loop/.style={min distance=8mm} %% Default loop has arrow.
}

\newtheorem{thm}{Theorem}[section]
\newtheorem{cor}[thm]{Corollary}
\newtheorem{lem}[thm]{Lemma}
\newtheorem{prop}[thm]{Proposition}

\newtheorem{obs}[thm]{Observation}

\theoremstyle{definition}
\newtheorem{rem}[thm]{Remark}

\theoremstyle{definition}
\newtheorem{defn}[thm]{Definition}

\theoremstyle{definition}
\newtheorem{ex}[thm]{Example}

\def\mtx#1{\begin{pmatrix} #1 \end{pmatrix}}

\DeclareMathOperator{\rank}{rank}
\DeclareMathOperator{\diag}{diag}
\DeclareMathOperator{\tr}{tr}
\DeclareMathOperator{\Span}{span}

\newcommand{\R}{\mathbb{R}}

\newcommand{\N}{\mathbb{N}}

\newcommand{\Rnn}{\R^{n\times n}}

\newcommand{\T}{\mathcal{T}}
\newcommand{\Mc}{\mathcal{M}}

\newcommand{\bq}{{\bf x}} % can't use q due to q = # evals
\newcommand{\br}{{\bf y}}
\newcommand{\bx}{{\bf x}}
\newcommand{\bv}{{\bf v}}

\newcommand{\bz}{{\bf z}}
\newcommand{\bu}{{\bf u}}
\newcommand{\bb}{{\bf b}}
\newcommand{\ba}{{\bf a}}
\newcommand{\bc}{{\bf c}}
\newcommand{\by}{{\bf y}}

\newcommand{\bzero}{{\bf 0}}
\newcommand{\oml}{{\bf m}}

\newcommand{\x}{\times}
\newcommand{\bit}{\begin{itemize}}
\newcommand{\eit}{\end{itemize}}
\newcommand{\ben}{\begin{enumerate}}
\newcommand{\een}{\end{enumerate}}
\newcommand{\beq}{\begin{equation}}
\newcommand{\eeq}{\end{equation}}
\newcommand{\bea}{\begin{eqnarray*}}
\newcommand{\eea}{\end{eqnarray*}}
\newcommand{\bean}{\begin{eqnarray}}
\newcommand{\eean}{\end{eqnarray}}
\newcommand{\bpf}{\begin{proof}}
\newcommand{\epf}{\end{proof}\ms}
\newcommand{\bmt}{\begin{bmatrix}}
\newcommand{\emt}{\end{bmatrix}}
\newcommand{\ms}{\medskip}

\newcommand{\beqa}{\begin{array}}
\newcommand{\eeqa}{\end{array}}

\newcommand{\lc}{\left\lceil}
\newcommand{\rc}{\right\rceil}
\newcommand{\lf}{\left\lfloor}
\newcommand{\rf}{\right\rfloor}
\newcommand{\wt}{\widetilde}

\newcommand{\mptn}{\mathcal{S}} %pattern manifold
\newcommand{\mrk}{\mathcal{R}} %rank manifold
\newcommand{\mei}{\mathcal{E}} %spectrum manifold
\newcommand{\mmu}{\mathcal{U}} %multiset pertibation manifold
\newcommand{\Rsn}{{S}_n(\R)} %symmetric matrix space 
\newcommand{\Rskewn}{{K}_n(\R)} %symmetric matrix space 
\newcommand{\G}{\mathcal{G}}

\newcommand{\vect}{\operatorname{vec}}

\DeclareMathOperator{\Mplus}{M_{+}}
\DeclareMathOperator{\mr}{mr}
\DeclareMathOperator{\M}{M}

\def \dunion {\dot{\bigcup}}
\def\spec{\operatorname{spec}}
\def\N{\mathcal{N}}
\def\man {\mathcal{M}}
\def\trans{^{\top}}
\newcommand{\qs}{q_S}
\newcommand{\qm}{q_M}
\newcommand{\Pc}{\mathcal P}
\newcommand{\Tc}{\mathcal T}

\begin{document}
\title{Generalizations of the Strong Arnold Property and  the minimum number of distinct eigenvalues of  a graph}
\author{Wayne Barrett\thanks{Department of Mathematics, Brigham Young University, Provo UT 84602, USA (wb@mathematics.byu.edu, H.Tracy@gmail.com).}
\and Shaun Fallat\thanks{Department of Mathematics and Statistics, University of Regina, Regina, Saskatchewan, S4S0A2, Canada (shaun.fallat@uregina.ca).}
 \and H.
Tracy Hall\footnotemark[1]
\and  Leslie
Hogben\thanks{Department of Mathematics, Iowa State University,
Ames, IA 50011, USA (LHogben@iastate.edu) and American Institute of Mathematics, 600 E. Brokaw Road, San Jose, CA 95112, USA
(hogben@aimath.org).}
\and Jephian C.-H. Lin\thanks{Department of Mathematics, Iowa State University,
Ames, IA 50011, USA (chlin@iastate.edu).}
\and Bryan L.  Shader\thanks{Department
of Mathematics, University of Wyoming,  Laramie, WY 82071, USA
(bshader@uwyo.edu).}  
 }

%\linenumbers
 
%%%%%%%%%%%%%%%%%%%%%%%%%%%%%%%%%%%%%%%%%%%%%

\maketitle

%%%%%%%%%%%%%%%%%%%%%%%%%%%%%%%%%%%%%%%%%%%%
\begin{abstract}
For a given graph $G$ and an associated class of real symmetric matrices whose off-diagonal entries are governed by
the adjacencies in $G$,  the collection of all possible spectra for such matrices is considered. Building on the pioneering
work of Colin de Verdi\`ere in connection with the Strong Arnold Property,  two extensions  are devised that target a better understanding of all possible spectra and their associated multiplicities. These new properties are referred to as the Strong Spectral Property and the Strong Multiplicity Property. Finally, these ideas are applied to  the minimum number of distinct eigenvalues associated with $G$, denoted by $q(G)$.  The  graphs for which $q(G)$ is at least  the number of vertices of $G$ less one are characterized.
\end{abstract}

\medskip


%%%%%%%%%%%%%%%%%%%%%%%%%%%%%%%%%%%%%%%%%%%%%
%%%%%%%%%%%%%%%%%%%%%%%%%%%%%%%%%%%%%%%%%%%%

\section{Introduction}\label{sintro}
Inverse eigenvalue problems appear in various contexts throughout mathematics and engineering.
The general form of an inverse eigenvalue problem is the following:\  given a family $\mathcal{F}$ of matrices and 
a spectral property $\Pc$, determine if there exists a matrix $A \in \mathcal{F}$ with property $\Pc$. Examples of families are tridiagonal, Toeplitz, or all symmetric matrices with a given graph.  Examples of properties include:\ having a prescribed rank, a prescribed spectrum,  a prescribed eigenvalue and corresponding eigenvector, or a prescribed 
list of multiplicities. Our focus is on  inverse eigenvalue problems where $\mathcal{F}$ is a set of symmetric matrices associated with a graph.  These have  received considerable attention, and a rich mathematical theory has been developed around them  (see, for example, \cite{FH}). 

All matrices in this paper are real.  Let  {$G=(V(G),E(G))$} be a (simple, undirected) graph  with  vertex set ${V(G)}=\{1,\dots,n\}$ and edge set {$E(G)$}.
The {\em set $\mptn(G)$ of symmetric matrices described by $G$} consists of the set of all  symmetric $n \times n$ matrices $A=(a_{ij})$  such that  for $i\ne j$, $a_{ij} \neq 0$ if and only if $ij \in {E(G)}$.  We denote  the spectrum of $A$, i.e.,  the multiset  
of eigenvalues of $A$, by $\spec(A)$.  The {\em inverse spectrum problem for G}, also known as the {\em inverse eigenvalue problem for G}, refers to determining the possible spectra that occur among the matrices in $\mptn(G)$. The inverse spectrum problem for $G$
seems to be difficult, as evidenced by that fact that it has been completely solved for only a few special families of graph, e.g. paths, generalized stars, double generalized stars, and complete graphs \cite{BF05, BLMNSSSY13, JLS03}. 

To gain a better understanding of the inverse eigenvalue problem for graphs, other spectral properties have been studied.
For example, the  {\it maximum multiplicity  problem for $G$}  is: Determine $\M(G)$, where
\[
\M(G) = \max\{\mbox{mult}_A (\lambda) :  A \in \mptn(G),\;\lambda\in \spec(A)\},
\]
and mult$_A(\lambda)$ denotes the multiplicity of $\lambda$ as an eigenvalue of $A$. A related  invariant is the  minimum rank of $G$, which is defined by 
 \[
 \mr(G) = \min\{\rank A :  A \in \mptn(G)\}.
 \]
 The {\it minimum rank problem} is:\ Given a graph $G$, determine $\mr(G)$.
As $\mr(G)+\M(G)=|G|$, where $|G|$ denotes the number of vertices of $G$, the maximum multiplicity and minimum rank problems
are essentially the same.  These problems have been extensively studied in recent years; see  \cite{FH, HLA2} for surveys. If the distinct eigenvalues of $A$ are $\lambda_1< \lambda_2< \cdots< \lambda_q$ and the multiplicity 
of these eigenvalues are $m_1, m_2, \ldots, m_q$ respectively, then the {\it ordered multiplicity list of $A$} 
is $\oml=(m_1, m_2, \ldots, m_q)$.  This notion gives rise to the {\em inverse ordered multiplicity list  problem}:\ Given a graph $G$, determine
which ordered multiplicity lists arise among the matrices in $\mptn(G)$.  This problem has been studied in \cite{BF04, BF05, JLS03}.
A recently introduced spectral problem  (see \cite{AACFMN13}) is the  {\em minimum number of distinct eigenvalues problem}:\ Given a graph $G$,  determine $q(G)$ where 
$q(G) = \min\{q(A)\; \mbox{:}\; A \in \mptn(G)\}$ and $q(A)$ is the number of distinct eigenvalues of $A$.

For a specific graph $G$  and a specific property $\Pc$, it is often difficult to find an explicit matrix $A \in \mptn (G)$ having property $\Pc$
(e.g., consider the challenge of finding a  matrix  whose graph is a path on five vertices and that has eigenvalues $0$, $1$, $2$, $3$, $4$).  

In this paper, we carefully describe the underlying theory of a technique based on the implicit function theorem, and develop new methods for two types of inverse eigenvalue problems. Suppose $G$ is a graph and $\Pc$ is a spectral property (such as having a given spectrum, given ordered multiplicity list, or given multiplicity of the eigenvalue 0) of a matrix in $\mptn(G)$. The  theory is applied to determine conditions (dependent on the property $\Pc$) that  guarantee if a matrix  $A\in\mptn(G)$ has property $\Pc$ and satisfies these conditions, then for every supergraph $\widetilde{G}$ 
  of $G$ (on the same vertex set) there is a matrix  $B\in \mptn(\widetilde{G})$ satisfying property $\Pc$. (A graph $\wt G$ is a {\em supergraph} of $G$ if $G$ is a subgraph of $\wt G$.)  In Section \ref{S2}, the technique is developed, and related to both the implicit function theorem and the work of 
Colin de Verdi\`ere \cite{CdV, CdV2}.  The technique is then applied to produce  two new properties for symmetric 
matrices called the Strong Spectral Property (SSP) and the Strong Multiplicity Property (SMP) that generalize the 
Strong Arnold Property (SAP). 

In Section \ref{commonSSPSMP}, we  establish general theorems for the inverse spectrum  problem, respectively the 
inverse multiplicity list problem, using matrices satisfying the  SSP, respectively the SMP. 
In Section \ref{sq}, we use the SSP and SMP to prove properties of $q(G)$.  In particular we  answer a question raised in  \cite{AACFMN13} by giving a complete characterization of the graphs 
$G$ for which $q(G)\geq |G|-1$.  %In Section 5, some results about $G$ for which $q(G)$ is small are derived.

%%%%%%%%%%%%%%%%%%%%%%%%%%%%%%%%%%%%%%%%%%%%
%%%%%%%%%%%%%%%%%%%%%%%%%%%%%%%%%%%%%%%%%%%%

\section{Motivation and fundamental results}\label{S2}
We begin by recalling an inverse problem due to Colin de Verdi\`ere.
 In his study of certain Schr\"odinger operators, Colin de Verdi\`ere was concerned with the maximum nullity, $\mu(G)$,
of matrices in the class of matrices $A \in \mptn (G)$ satisfying: 
\begin{itemize}
\item[(i)] all off-diagonal entries of $A$ are non-positive;
\item[(ii)] $A$ has  a unique negative eigenvalue with multiplicity $1$; and 
\item[(iii)] $O$ is the only symmetric matrix $X$ satisfying $AX=O$, $A\circ X=O$ and $I\circ X=O$.
\end{itemize}
Here $\circ$ denotes the Schur (also known as the Hadamard or entrywise) product and $O$ denotes the zero matrix. Condition (iii) is known as the {\it Strong Arnold Property}  (or SAP for short).
Additional variants, called {\em Colin de Verdi\`ere type parameters}, include $\xi(G)$, which is the maximum nullity of matrices $A \in \mptn (G)$ satisfying the SAP \cite{BFH05}, and $\nu(G)$, which is the maximum nullity of positive semidefinite matrices $A \in \mptn (G)$ satisfying the SAP \cite{CdV2}.  
The maximum multiplicity $\M(G)$ is also the maximum nullity of a matrix in $\mptn(G)$, so $\xi(G)\le \M(G)$. Analogously, $\nu(G)\le \Mplus(G)$, where $\Mplus(G)$  is the maximum nullity among positive semidefinite matrices in $\mptn(G)$. 

Colin de Verdi\`ere \cite{CdV} used results from manifold theory to show conditions equivalent to (i)--(iii) imply that 
%if $G$ is a minor\footnote{A minor of a graph is obtained by contracting edges and by deleting edges and vertices.} of $\widetilde{G}$, then there exists $B \in \mptn (\widetilde{G})$ such that $A$ and $B$ have the same nullity.  In particular, 
if $G$ is a subgraph of $\widetilde{G}$, then the existence of $A\in \mptn (G)$ satisfying (i)--(iii), 
implies the existence of $\widetilde{A} \in \mptn (\widetilde{G})$ having the same nullity as $A$ and satisfying (i)--(iii).\footnote{In fact, much more was shown: If $G$ is a minor of $H$, then there exists $B \in \mptn (H)$  satisfying (i)--(iii) such that $A$ and $B$ have the same nullity.  A {\em minor} of a graph is obtained by contracting edges and by deleting edges and vertices, and a graph parameter $\beta$ is {\em minor monotone} if $\beta(G)\le\beta(H)$.}  The formulation  of the SAP in (iii) is due to van der Holst, Lov{\'a}sz and Schrijver \cite{HLS}, which gives  a linear algebraic treatment of the SAP and the Colin de Verdi\`ere number $\mu$. 
Using a  technique similar to that in \cite{HLS}, Barioli, Fallat, and Hogben \cite{BFH05} 
showed that  if there exists $A\in \mptn (G)$ satisfying the SAP and  $G$ is a subgraph of  $\widetilde{G}$,  then there exists  $\widetilde{A} \in \mptn (\widetilde{G})$ such that $\wt A$ has the SAP and $A$ and $\widetilde{A}$ 
have the same nullity.\footnote{Again, the result was established for minors: If $G$ is a minor of $H$, then there exists $B \in \mptn (H)$  satisfying the SAP such that $A$ and $B$ have the same nullity.}  

The Strong Arnold Property has been used to obtain many results about the maximum nullity of a graph. %; for a linear algebraic treatment of SAP and the Colin de Verdi\`ere number $\mu$ see \cite{HLS}, and for an application of SAP to maximum nullity see \cite{BFH05}.  
Our goal in this section is to first describe the general technique behind Colin de Verdi\`ere's work, and then develop analogs of the SAP for the inverse 
spectrum problem and the inverse multiplicity list problem.
For   convenience, we state below a version of the Implicit Function Theorem (IFT) because it is central to the technique; see  \cite{DR, KP}.

\begin{thm}\label{IFT}
Let $F{\; :\;}\mathbb{R}^{s+r} \rightarrow \mathbb{R}^t$ be a continuously differentiable
 function on an open subset $U$ of 
 $\mathbb{R}^{s+r}$
defined by
\[
F( \bx, \by)=(F_1( \bx, \by), F_2( \bx, \by), \ldots, F_t( \bx, \by)),
\]
 where $ \bx=(x_1, \ldots, x_s)\trans \in \mathbb{R}^s$ and $ \by \in \mathbb{R}^r$. 
 Let $( \ba, \bb)$ be an element of $U$
 with $ \ba\in \mathbb{R}^s$ and $ \bb\in \mathbb{R}^r$,  and $ \bc\in \mathbb{R}^t$ such that 
  $F( \ba, \bb)= \bc$.
  If  the $t\times s$ matrix 
  \begin{equation}
  \label{Jac}
   \left(\frac{\partial F_i}{\partial x_{j}}\;{\rule[-3.6mm]{.1mm} {8mm}}_{\;( \ba, \bb)} \right) 
   \end{equation}
  has rank $t$, then there exist an open neighborhood $V$  containing $ \ba$
 and an open neighborhood $W$ containing $ \bb$ 
 such that  $V \times W \subseteq U$ and a continuous function $\phi: W \rightarrow V$
 such that 
$F( \phi(\by), \by)= \bc$ for all $\by \in W$. \end{thm}

The IFT concerns robustness of solutions to the system $F(\bx,\by)=\bc$.
Namely, the existence of a ``nice'' solution $\bx=\ba$ to $F(\bx, \bb)=\bc$ guarantees 
the existence of solutions to all systems $F(\bx,\tilde{\bb})=\bc$ with $\tilde{\bb}$ sufficiently close  to $\bb$.
Here, ``nice'' means that the  columns of   the matrix in (\ref{Jac}) span $\R^t$. 

%Insert right after Theorem 2.1
\begin{rem}\label{rem:IFT}
More quantitative proofs of the implicit function theorem (see \cite[Theorem 3.4.10]{KP}) show that there exists an $\epsilon> 0$ that depends only on $U$ and the Jacobian in (1) such that there is such a continuous function $\phi$ with $F(\phi(\by),\by)=\bc$
for all $\by$ with $\|\by-\bb\| <\epsilon$.
\end{rem}
   
A useful application of the IFT arises in the setting of smooth manifolds.  We refer the reader to \cite{KP,L} for basic definitions and results about manifolds.  Given a manifold $\man$ embedded smoothly in an inner product space and a point $\bx$ in $\man$, we denote the tangent space to $\man$ at $\bx$ by $\T_{\man.\bx}$ and the normal space\footnote{That is, the orthogonal complement of the tangent space in the ambient inner product space.}
to $\man$ at $\bx$ by $\N_{\man.\bx}$. For the manifolds of matrices discussed here, the inner product of $n\x n$ matrices $A$ and $B$ is $\langle A,B\rangle=\tr(A\trans B)$, or equivalently, the Euclidean inner product on $\R^{n^2}$ (with matrices viewed as $n^2$-tuples). 

Many problems, including our inverse problems, can be reduced to determining whether or not the intersection of two manifolds $\man_1$ and $\man_2$ is non-empty.  There is a condition known as transversality such that if $\man_1$ and $\man_2$ 
intersect transversally  at $\bx$, then any manifold ``near'' $\man_1$ and  any manifold ``near'' $\man_2$ intersect non-trivially.   In other words, 
the existence of a nice solution implies the existence of a solution to 
all ``nearby'' problems.   More precisely, let $\man_1$ and $\man_2$  be manifolds in some $\R^d$, and $\bx$  be a point in $\man_1 \cap \man_2$.  The manifolds $\man_1$ and $\man_2$ intersect  {\it transversally}  at $\bx$ provided $\T_{\man_1.\bx} + \T_{\man_2.\bx} = \R^d$, or equivalently
$\N_{\man_1.\bx} \cap \N_{\man_2.\bx}= \{O\}$.

By a smooth family $\man(s)$  $(s\in (-1,1))$ of manifolds in $\R^d$  we mean that 
$\man(s)$ is a manifold in $\R^d$  for each $s\in (-1,1)$, 
and $\man(s)$ varies smoothly as a function of $s$.   Thus, if $\man(s)$ is a smooth family of 
manifolds in $\R^d$, then

\begin{itemize}
\item[(i)] there is a $k$ such that $\man(s)$ has dimension $k$ for all $s\in (-1,1)$;
\item[(ii)] there exists a collection $\{U_{\alpha}: \alpha \in \cal{I} \}$ of relatively open sets  
whose union is $\cup_{s\in (-1,1) } \man(s)$;  and 
\item[(iii)] for each $\alpha$ there is a diffeomorphism $F_{\alpha}: (-1,1)^{k+1} \rightarrow U_\alpha$.
\end{itemize}

\noindent
Note that  if $\bx_0\in (-1,1)^k$, $s_0 \in (-1,1)$, $\by_0 \in \man(t_0)$ and $F_{\alpha}(\bx_0,s_0)=\by_0$, 
then the tangent space to $\man({s_0})$ at $y_0$ is the  column space of  the $d\times k$ matrix
\[
  \left(\frac{\partial F_{\alpha}}{\partial x_{j}}\;{\rule[-3.6mm]{.1mm} {8mm}}_{\;( \bx=\bx_0, s=t_0)} \right)\!.
\]

%it is a continuous function of $\epsilon$. 
The following theorem can be viewed as a  specialization of  Lemma 2.1 and Corollary 2.2   of \cite{HLS}  to the case of two manifolds.

%In loose terms the following theorem 
%says that if $M(0)$ and $N(0)$ intersect transversally at a point $x(0)$, then for all $\epsilon$ sufficiently small 
%$M(\epsilon)$ and $N(\epsilon)$ intersect and the point of intersection $x(\epsilon)$ can be chosen so that 

\begin{thm}
\label{thm:transverse}
 Let $\man_1(s)$ and $\man_2(t)$ be  smooth families  of manifolds in $\R^d$,  and assume that $\man_1(0)$ and $\man_2(0)$ intersect transversally at $\by_0$.  Then there is a neighborhood $W\subseteq \R^2$ of the origin and a continuous function $f:W \rightarrow \R^d$
  such that for each $\epsilon=(\epsilon_1, \epsilon_2) \in W$, $ \man_1(\epsilon_1) $ and $ \man_2(\epsilon_2)$ intersect transversally at $f(\epsilon)$.
\end{thm} 

\bpf
Let $k$, respectively $\ell$, be the common dimension of each of the  manifolds $\man_1(s)$ $(s\in (-1,1))$ and  
each of the manifolds $\man_2(t)$ $(t\in (-1,1))$, respectively. 
Let $\{ U_{\alpha}: \alpha \in \mathcal{I} \}$  and $ F_\alpha:(-1,1)^{k+1} \rightarrow U_\alpha $  be the open sets and diffeomorphisms
for $\man_1(t)$ discussed above.
Let $\{ V_{\beta}: \beta \in \mathcal{J} \}$ and  $G_{\beta}: (-1,1)^{\ell+1} \rightarrow V_\beta$
 be the similar open sets and diffeomorphisms for $\man_2(t)$.
 
Choose $\alpha$ so that $\by_0\in U_{\alpha}$ %contains $\man_1(0)$ 
and $\beta$ so that $\by_0\in V_{\beta}$. 
 Define 
\[
H: \R^k \times (-1,1) \times \R^{\ell} \times (-1,1) \rightarrow \R^d \mbox{ by }  H(\bu,s,\bv,t)= F_\alpha(\bu,s) -G_\beta(\bv,s),
\]
where $\bu\in (-1,1)^k$, $\bv\in (-1,1)^{\ell}$, and $s,t \in (-1,1)$.
There exists $\bu(0)$ and $\bv(0)$ such that 
$F_\alpha(\bu(0),0)= \by_0$ and $G_\beta(\bv(0),0)=\by_0$.
Hence $H(\bu(0),0,\bv(0),0)=0$.
Since $F$ and $G$ are diffeomorphisms, the Jacobian of the function $H$ restricted to $s=0$ and $t=0$ 
and evaluated at $\bu=\bu(0)$ and $\bv=\bv(0)$ is the $d$ by $k+ \ell$ matrix 
\[
\mbox{Jac} =\left(  \begin{array}{cc|cc}
\phantom{Z} & \\
& \frac{\partial F_i}{\partial u_{j}}\;{\rule[-3.6mm]{.1mm} {8mm}}_{\;( \bu(0),0)}  & 
\phantom{Z}&\frac{\partial G_i}{\partial v_{j}}\;{\rule[-3.6mm]{.1mm} {8mm}}_{\;(\bv(0),0)}  \\
 \phantom{Z} &
 \end{array} \right)\!.
\]
The assumption that $\man_1(0)$ and $\man_2(0)$ intersect transversally at $\by(0)$ 
implies that the column space of $\mbox{Jac}$ is all of $\R^d$.  

The result now follows by applying the implicit function theorem (Theorem \ref{IFT}) and the fact that every 
matrix in a sufficiently small neighborhood of a full rank matrix has full rank. %  of a spanning set is a spanning set.
\epf

%%%%%%%%%%%%%%%%%%%%%%%%%%%%%%%%%%%%%%%%%%%%

\subsection{The Strong Arnold Property}

We now use Theorem \ref{thm:transverse} to describe a  tool for the inverse rank problem. This 
tool is described in \cite{HLS} and used there and  in \cite{param, GCC, BFH05}. We use this problem
to familiarize the reader with the general technique.

Let $G$ be a graph of order $n$, and $A\in \mptn(G)$. 
For this problem the two manifolds that concern us 
are $\mptn (G)$ and  the manifold
\[\mrk_A=\{B\in\Rsn:\rank(B)=\rank(A)\},\]
consisting of all $n\times n$ symmetric matrices with the same rank as $A$.
We view both of these as subsets of $\Rsn$, the set of all $n\times n$ symmetric 
matrices endowed with the inner product $\langle V, W \rangle = \tr(VW)$.
Thus $\mptn(G)$ and $\mrk_A$ can be thought of as submanifolds of $\R^{n(n+1)/2}$.
It is easy to see that 
\begin{eqnarray*}
\T_{\mptn({G}).A}&=& \{ X \in \Rsn: 
x_{ij} \neq 0 \mbox{$ \implies$ $ij$ is an edge of $G$ or $i=j$ } \} , \mbox{ and } \\
\N_{{\mptn}{(G).A}}& =&  \{ X \in \Rsn: A\circ X=O  \mbox{ and }  I\circ X=O \}.
\end{eqnarray*}

For $\N_{\mrk_{A}.A}$ we have the following result \cite{HLS}.
\begin{lem}
\label{lem:SAP}
Let $A$ be a symmetric $n\times n$ matrix of rank $r$. 
Then 
\[
\N_{\mrk_{A}.A}=\{ X \in \Rsn: AX=O\}.
\]
\end{lem}
\bpf
There exists an invertible $r\times r$ principal submatrix  of $A$, and 
by permutation similarity we may take this to be the leading principal submatrix. Hence $A$ has the form 
\[
\renewcommand*{\arraystretch}{1.5}
\left( \begin{array}{c|c}
A_1 & A_1U
\\ [.1pt] \hline 
U\trans A_1 & U\trans A_1 U 
\end{array} 
\right)
\]
for some invertible $r\times r$ matrix $A_1$ and some $r \times (n-r)$ matrix $U$.

Let $B(t)$ $(t\in (-1,1))$ be a differentiable path of symmetric rank $r$ matrices 
such that $B(0)=A$.  For $t$ sufficiently small, the leading $r\times r$ principal submatrix 
of $B(t)$ is invertible and $B(t)$ has the  form
\[
\renewcommand*{\arraystretch}{1.5}
\left( \begin{array}{c|c}
B_1(t) & B_1(t)U(t)\\ \hline 
U(t)\trans B_1(t) & U(t)\trans B_1(t) U(t) 
\end{array} 
\right)\!,
\]
where $B_1(t)$ and $U(t)$ are differentiable, $B_1(0)=A_1$ and $U(0)=U$.
Differentiating with respect to $t$ and then evaluating at $t=0$ gives 
\[
\renewcommand*{\arraystretch}{1.5}
\dot{B}(0)= 
\left( \begin{array}{c|c} 
\dot{B}_1(0) & \dot{B}_1(0)U \\ \hline 
U{\trans} \dot{B}_1(0) & U{\trans} \dot{B}_1(0) U
\end{array} \right)
+ 
\left( \begin{array}{c|c}
O & A_1 \dot{U}(0) \\ \hline
\dot{U}(0){\trans}A_1 &  \dot{U}(0){\trans} A_1U + U{\trans} A_1 \dot{U} (0)
\end{array} \right)\!.
\]
It follows that $\T_{\mrk_A.A}=\Tc_1+\Tc_2$, where
\[
\renewcommand*{\arraystretch}{1.5}
\Tc_1:=\left\{
\left( \begin{array}{c|c} 
R & RU \\ \hline 
U{\trans} R & U{\trans} R U 
\end{array} \right): R \mbox{ is an arbitrary symmetric  $r\x r$ matrix}  \right\}
\]
and
\[
\renewcommand*{\arraystretch}{1.5}
\Tc_2:=\left\{ \left( \begin{array}{c|c}
O & A_1S \\ \hline 
S{\trans}A_1 &  S{\trans} A_1U + U{\trans} A_1 S 
\end{array} \right): S \mbox{ is an arbitrary $r \times (n-r)$ matrix} \right\}.
\]
Consider an $n\times n$ symmetric matrix
\[
\renewcommand*{\arraystretch}{1.5}
W= \left( \begin{array}{cc} 
C & D \\ %\hline
D{\trans} & E 
\end{array}
\right)\!,
\]
where $C$ is $r \times r$. 
Then $W \in \Tc_2^{\perp}$ 
if and only if 
\[
\tr ( DS{\trans}A_1+ D{\trans}A_1S + ES{\trans}A_1U+ EU{\trans}A_1S) =0 \mbox{ for all $S$}
\]
or equivalently
\[
\tr ((D{\trans}A_1+ EU{\trans}A_1)S)=0 \mbox{ for all $S$.} 
\]
Thus, $W\in \Tc_2^{\perp}$ if and only if $D{\trans}A_1+EU{\trans}A_1=O$, which is equivalent to $D{\trans}=-EU{\trans}$ since $A_1$ is invertible.
Similarly, $W\in \Tc_1^{\perp}$ if and only if $C+DU\trans +UD\trans +UEU\trans$ is skew symmetric. As $C+DU\trans +UD\trans +UEU\trans$ is symmetric, 
we have $C\in \Tc_1^{\perp}$ if and only if $C+DU\trans +UD\trans +UEU\trans=0$.  For $W\in \Tc_1^{\perp}\cap \Tc_2^{\perp}$, $C=UEU\trans$ and therefore
\beq\label{neweq31a}\small
\renewcommand*{\arraystretch}{1.5}
\N_{\mrk_A.A}=\left\{ \left( \begin{array}{c|c}
UEU{\trans} & -UE \\ \hline
-EU{\trans} & E  \end{array} \right): \mbox{ 
 $E$ is an arbitrary symmetric $(n-r) \times (n-r)$ matrix} \right\}.
\eeq
It is easy to verify that this is precisely that set of symmetric matrices $X$ such that $AX=O$.
\epf

Lemma \ref{lem:SAP} implies that  $\mrk_A$ and $\mptn (G)$ intersect transversally  at $A$ 
if and only if $A$ satisfies the SAP. Now that we know the pertinent tangent spaces for the inverse rank problem for $G$,
we can apply Theorem \ref{thm:transverse} to easily obtain the following useful tool, which was previously 
proved and used in  \cite{HLS}.

\begin{thm}
\label{thm:SAP}
If $A\in \mptn(G)$ has the SAP, then every supergraph of $G$ with the same vertex set has a realization that has the same rank as $A$ and has the SAP.  
\end{thm}
\bpf
%Let $M$ be the manifold $\mrk_A$.  
Let $\wt G$ be a supergraph of $G$ with the same vertex set, and define $\Delta=E(\wt G)-E(G)$.  For $t\in (-1,1)$, define the manifold $\man(t)$ as those $B=(b_{ij})\in\Rsn$ with $b_{ij}\neq 0$ if $ij\in E(G)$, $b_{ij}=t$ if $ij\in \Delta$, and $b_{ij}=0$ if  $i\ne j$ and $ij\not\in E(\wt G)$.  Since $\man(0)=\mptn(G)$, and $A$ has the SAP,   $\mrk_A$ and $\man(0)$ intersect transversally at $A$.  Therefore Theorem  \ref{thm:transverse}
guarantees a continuous function $f$ such that for $\epsilon$ sufficiently small,
 $\mrk_A$ and $\man(\epsilon)$ intersect transversally at  $f(\epsilon)$, so   $f(\epsilon)$ has the same rank as $A$ and $f(\epsilon)$ has the SAP.   For $\epsilon>0$,  
$f(\epsilon) \in \mptn(\wt G)$.  
%Furthermore, $A_\epsilon$ also has SSP.
\epf

%%%%%%%%%%%%%%%%%%%%%%%%%%%%%%%%%%%%%%%%%%%%

\subsection{The Strong Spectral Property}

We now follow the general method outlined in the previous subsection to derive
an analog of the SAP and Theorem \ref{thm:SAP} for the inverse spectrum problem.  Certain aspects of this section were in part motivated by discussions with Dr. Francesco Barioli in connection with the inverse eigenvalue problem for graphs \cite{Bar}.

 Given a multiset $\Lambda$ of real numbers that has cardinality $n$, 
 the set of all $n\times n$ symmetric matrices with spectrum $\Lambda$ is denoted by 
$
\mei_{\Lambda}.
$
Thus  if $A\in \mei_{\Lambda}$, then 
$\mei_{\Lambda}$ is all symmetric matrices cospectral with $A$. 
It is well known that $\mei_{\Lambda}$ is a manifold \cite{Arnold}. 
A comment on notation: The notation $\mei_{\Lambda}$ for the constant spectrum manifold was chosen because this manifold is determined  by ${\Lambda}$.  Then a symmetric matrix $A$ is in $\mei_{\spec(A)}$.  In conformity with this, the constant rank manifold containing $A$ should be denoted $\mrk_{\rank A}$, but we follow the literature in denoting it by $\mrk_A$.

Let $A$ be a symmetric $n\times n$ matrix. The centralizer of $A$ is the set of all matrices that commute with $A$, and is denoted by 
$\mathcal{C}(A)$.   The commutator, $AB-BA$ of two matrices is denoted by $[A,B]$.   
 The next result is well known but we include a brief proof for completeness.
 \begin{prop}\label{centralizer} Suppose $A\in\Rsn$ and $\bv_i, i=1,\dots,n$ is an orthonormal basis of eigenvectors for $A$ with $A\bv_i=\mu_i\bv_i$ (where the $\mu_i$ need not be distinct). Then
  \[
  \mathcal{C}(A) =\Span (\{ \bv_i\bv_j{\trans} : \mu_i=\mu_j\}).
  \]
  \end{prop}
 \bpf Let $S=\Span (\{ \bv_i\bv_j{\trans} : \mu_i=\mu_j\})$. Clearly $S\subseteq \mathcal{C}(A)$.  For the reverse inclusion, observe that $\{ \bv_i\bv_j{\trans} : i=1,\dots,n, j=1,\dots,n\}$ is a basis for $\Rnn$.  Thus any $B\in \mathcal{C}(A)$ can be expressed as $B=\sum_{i=1}^n\sum_{j=1}^n \beta_{ij}\bv_i\bv_j\trans$. Then 
 \[O=[A,B]=\sum_{i=1}^n\sum_{j=1}^n\beta_{ij}(\mu_i-\mu_j)\bv_i\bv_j\trans.\]
By the independence of the $\bv_i\bv_j\trans$, $\beta_{ij}(\mu_i-\mu_j)=0$ for all $i,j$.  Therefore, $\mu_i\ne\mu_j$ implies  $\beta_{ij}=0$, so $\mathcal{C}(A)\subseteq S$.  \epf
 
Throughout this section we assume that   $\lambda_1, \ldots, \lambda_q$ are the distinct eigenvalues of $A$
and that  $A=\sum_{i=1}^q \lambda_iE_i$ is the spectral decomposition of $A$ (i.e.~the $E_i$ are 
mutually orthogonal idempotents that sum to $I$). %; if $\bv$ is an eigenvector of $A$ for eigenvalue $\lambda_i$, then $E_i\bv=\bv$ and $E_j\bv=0$. 
\begin{lem}
\label{lem:SSP}
Let $A$ be a symmetric $n\times n$ matrix with $\mbox{spec}(A)=\Lambda$.  Then
$\N_{\mei_{\Lambda}.A}= \mathcal{C}(A) \cap \Rsn$.
\end{lem}

\bpf
Consider a differentiable path $B(t)$ $( t\in (-1,1))$ on $\mei_{\Lambda}$ such that $B(0)= A$.
Then $B(t)$ has spectral decomposition
\[
B(t) = \sum_{i=1}^{q} \lambda_i F_i (t).
\]
Clearly  $F_i(0)=E_i$ for $i=1,2,\ldots, q$.  As $F_i(t)$ is given by a polynomial in $B(t)$ with coefficients that depend only on the spectrum of $B(t)$, i.e. on $\Lambda$,
\cite[Corollary to Theorem 9 (Chapter 9)]{HK}, 
 $F_i(t)$ is a differentiable function of $t$.
Since  the  $F_i(t)$ are mutually orthogonal idempotents, we have that 
\begin{eqnarray}
\label{eq:idem1}
\dot{F}_i(0) E_i+ E_i \dot{F}_i(0) & = & \dot{F}_i(0) \mbox{\ \ \ and } \\
\label{eq:idem2}
\dot{F}_i(0)E_j+ E_i \dot{F}_j (0)&=&O \qquad \mbox{ $i,j=1,2,\ldots, q$ and $i\neq j$.}
\end{eqnarray}
Post-multiplying (\ref{eq:idem1}) by $E_j$ 
gives
\begin{equation}
\label{obs1}
E_i\dot{F}_i(0)E_j  = \dot{F}_i(0)E_j \mbox{ \qquad \ \, for all $j\neq i$}
\end{equation}
Pre-multiplying (\ref{eq:idem2}) by $E_j$ gives
\begin{equation}
\label{obs2}
E_j\dot{F}_i(0)E_j  = O \qquad \qquad \quad \mbox{\;\,for all $j\neq i$.}
\end{equation} 
Post-multiplying (\ref{eq:idem1}) by $E_i$ gives 
\begin{equation}
\label{obs3}
E_i\dot{F}_i(0)E_i=O \qquad \qquad \qquad \mbox{ \; \; for all $i$}.
\end{equation}
Equation (\ref{obs1}) %and the fact that $E_j$ is idempotent 
implies that  the image of $\dot{F}_i(0)E_j$ is contained in the image of $E_i$ $(j\neq i)$.

 Consider eigenvectors $\bq$ and $\br$ of $A$.  First suppose that 
they  correspond to the same eigenvalue, say $\lambda_i$.
Since $\bq=E_i\bq$ and $\br=E_i\br$, \eqref{obs2} and (\ref{obs3}) imply that $\br{\trans} \dot{F}_j(0) \bq=0$ for all $j$. 
Thus
\begin{eqnarray*}
\tr\left(\dot{B}(0)(\bq\br{\trans} + \br\bq{\trans})\right) & = & \tr\left(\sum_{j=1}^q \lambda_j \dot{F}_j(0)(\bq\br{\trans} + \br\bq{\trans})\right) \\
& = & \sum_{j=1}^q \lambda_j 2\br{\trans} \dot{F}_j(0) \bq\\
& = & 0.
\end{eqnarray*}
Therefore $\bq\br{\trans} +\br\bq{\trans} \in \N_{\mei_{\Lambda}.A} $ for all such choices of $\bq$ and $\br$ corresponding to the same eigenvalue.

Now suppose that $\bx$ and $\by$ are unit eigenvectors of $A$ corresponding to distinct eigenvalues, say $\lambda_1$ and $\lambda_2$.  Let $\mu_3 , \dots , \mu_n$ be the remaining eigenvalues of $A$ (with no assumption that they are distinct from each other or $\lambda_1$ or $\lambda_2$) and $\bz_3, \dots , \bz_n$ a corresponding orthonormal set of eigenvectors.  Let 
\[
D(t) =\left(\begin{array}{rr} \cos t & \sin t \\ -\sin t & \cos t \end{array}\right) \left(\begin{array}{rr} \lambda_1 & 0 \\ 0 & \lambda_2  \end{array}\right) \left(\begin{array}{rr} \cos t & -\sin t \\ \sin t & \cos t \end{array}\right) \oplus \diag \, (\mu_3 , \dots , \mu_n),
\]
let $S=\left(\begin{array}{ccccc}  \bx & \by & \bz_3 & \cdots & \bz_n \end{array}\right)$, and let $B(t) = S D(t) S{\trans}$.  Then $B(t)$ is cospectral with $A$, $B(0)= A$, and 
\[
\dot{B}(0)=S \dot{D}(0) S{\trans} =\begin{pmatrix} \bx & \by \end{pmatrix}  \begin{pmatrix}0 & \lambda_2 - \lambda_1\\ \lambda_2 - \lambda_1 & 0\end{pmatrix} \begin{pmatrix} \bx{\trans} \\ \by{\trans} \end{pmatrix} =  (\lambda_2 - \lambda_1) ( \bx \by{\trans} + \by \bx{\trans} ).
\]
It follows that $\bx\by{\trans}+\by\bx{\trans}\in \T_{\mei_{\Lambda}.A}$  for any eigenvectors $\bx$ and $\by$ corresponding to distinct eigenvalues of $A$.

Let $\bv_1, \ldots, \bv_n$ be a basis of eigenvectors of $A$. Then 
$\{\bv_i\bv_j{\trans} + \bv_j\bv_i{\trans}: i \leq j \}$ forms a basis for $\Rsn$.
We have shown that if $\bv_i$ and $\bv_j$ correspond to distinct eigenvalues, then
 $\bv_i\bv_j{\trans} + \bv_j\bv_i{\trans}\in \T_{\mei_{\Lambda}.A}$   and if $\bv_i$ and $\bv_j$ correspond to the same eigenvalues, then 
  $\bv_i\bv_j{\trans} + \bv_j\bv_i{\trans}\in \N_{\mei_\Lambda.A}$.  
  It follows that:
  \bean 
   \T_{\mei_{\Lambda}.A}\!\!&=\!\!&\Span (\{\bv_i\bv_j{\trans} + \bv_j\bv_i{\trans}: \mbox{$\bv_i$ and $\bv_j$ correspond to distinct eigenvalues of } A \}).  \label{neweq31z}\\
  \N_{\mei_{\Lambda}.A}\!\!&=\!\!&\Span (\{\bv_i\bv_j{\trans} + \bv_j\bv_i{\trans}: \mbox{$\bv_i$ and $\bv_j$ correspond to the same eigenvalue of } A \}).  \label{neweq31b}\eean
 By Proposition \ref{centralizer}, $\mathcal{C}(A)$ is the span of the set of matrices of the form $ \bv_i\bv_j{\trans}$ where $\bv_i$ and $\bv_j $ correspond to the same eigenvalue of   $A$.
  Thus $ \N_{\mei_{\Lambda}.A}=\mathcal{C}(A) \cap \Rsn$.
  \epf

\begin{defn} The symmetric matrix $A$ has the {\it Strong Spectral Property} (or $A$ has the SSP for short)
if the only symmetric matrix $X$ satisfying $A\circ X=O$, $I\circ X=O$ and $[A,X]=O$ is $X=O$. 
\end{defn}

\begin{obs} For symmetric matrices $A$ and $X$, $[A,X]=O$ if and only if $AX$ is symmetric, so the SSP could have been defined as: The only symmetric matrix $X$ satisfying $A\circ X=O$, $I\circ X=O$ and $AX$ is symmetric is $X=O$.  
\end{obs}
Lemma \ref{lem:SSP}
asserts that $A$ has the SSP if and only if  the manifolds $\mptn(G)$ and $\mei_{\Lambda}$ intersect transversally at $A$, where $G$ is the graph such that $A\in\mptn(G)$ and $\Lambda=\spec(A)$. 
  A proof similar to that of
Theorem \ref{thm:SAP} yields the next result.

\begin{thm}
\label{thm:SSP}
If $A\in \mptn(G)$ has the SSP, then every supergraph of $G$ with the same vertex set has a realization that has the same spectrum as $A$ and has the SSP.
\end{thm}

For every $A\in\mptn(K_n)$, $A\circ X=O$ and $I\circ X=O$ imply $X =O$, so trivially $A$ has the SSP. Next we discuss some additional examples.  

\begin{ex}
\label{exstar}
Let
\[
A=\left(\begin{array}{cccc}0&1&1&1\\1&0&0&0\\1&0&0&0\\1&0&0&0\\ \end{array} \right)
.
\]
Then $\spec(A)=\{0,0,\sqrt{3}, -\sqrt{3}\}$.  To show that $A$ has the SSP, 
consider $X\in S_4(\R)$ such that $A\circ X=0$, $I \circ X=O$ and $[A,X]=O$.   
Then  $X$   is a matrix of the form 
\[
X=\left( \begin{array}{ccccc} 0&0&0&0\\0&0&u&v\\0&u&0&w\\0&v&w&0\end{array} \right). 
\]
The fact that $X$ commutes with $A$ implies that $X$ has all row sums and column sums equal to zero, which in turn 
implies $X=O$.  Thus, $A$ has the SSP, and  by Theorem \ref{thm:SSP}, every supergraph of $K_{1,3}$ has a realization with spectrum $\{0,0,\sqrt{3}, -\sqrt{3}\}$.
\end{ex}

\begin{ex}
Let $G$ be the star on $n\ge 5$ vertices having $1$ as a central vertex, 
 let $A\in \mptn(G)$, and let $\lambda$ be an eigenvalue of $A$ of multiplicity at least 3.
From a theorem of Parter and Wiener (see, for example, \cite[Theorem 2.1]{FH}),  $\lambda$ occurs on the diagonal of $A(1)$  at least 4 times.\footnote{$A(1)$ is the principal submatrix of $A$ obtained by deleting row and column 1.}
Without loss of generality, we may assume that $a_{22}= a_{33}=a_{44}=a_{55}= \lambda$.
Let $\bu= [0 ,0 ,0,a_{15},-a_{14},0, \dots ,0]\trans$\!, $\bv= [0 , a_{13} ,-a_{12},0 ,0,0, \dots ,0]\trans$\!,
%\[u= \left( \begin{array}{c}0 \\0 \\0\\a_{15}\\-a_{14} \\0 \\\vdots \\0\end{array}\right),  \quad 
%v= \left(\begin{array}{c}0 \\ a_{13} \\-a_{12} \\0 \\0 \\0 \\\vdots \\0 \end{array} \right), \]
and $X=\bu{\bv}\trans+\bv{\bu}\trans$\!.
It can be verified that $AX= \lambda X= XA$, $A\circ X=O$ and $I\circ X=O$.
Thus, $A$ does not have the SSP.
Therefore, no matrix in $\mptn(G)$ with an eigenvalue of multiplicity at least $3$
has the SSP.
\end{ex} 

\begin{ex}
Let $G$ be as in the previous example.
Let $A\in \mptn(G)$ such that no eigenvalue of $A$ has multiplicity $3$ or more.
Without loss of generality we may assume that 
$A(1)= \oplus_{j=1}^k \lambda_j I_{n_j}$
for some distinct $\lambda_1, \dots, \lambda_k$ and positive integers $n_1, \dots, n_k$ 
with $n-1=n_1+ n_2 + \cdots +n_k$.  As every eigenvalue of $A$ has multiplicity $2$ or less,
each $n_j \leq 3$ for $j=1,2,\ldots, k$. 

Let $X$  be a symmetric matrix with $[A,X]=O$,  $A\circ X=O$ and $I\circ X=O$.
The last two conditions imply that   all entries in the first row, the first column, and on the diagonal of $X$ are zero.
This and the first condition imply that $X(1)$ is in $\mathcal{C}(A(1))$.  By the distinctness of the $\lambda_j$
we conclude that $X(1)= \oplus_{j=1}^k X_j$ where $X_j$ is a symmetric matrix of order $n_j$ with zeros on 
the diagonal. Partition $A$ as $A=\mtx{\alpha & \ba\trans\\ \ba & A(1)}$ and the partition $\ba=(\ba_1\trans,\dots,\ba_k\trans)\trans$ conformally with $X(1)$.  Then $[A,X]=O$ implies that  $X_j\ba_j=\bzero$.  As $n_j \leq 3$, $X_j$ is symmetric
and has zeros on its diagonal and every entry of $\ba_j$ is nonzero, this implies $X_j=O$. Thus, $X=O$ and we conclude that $A$ has the SSP.
\end{ex}
%\begin{ex}
%Let $A\in \mptn(P_n)$.  
%Suppose $X$ is a symmetric matrix such that $X\circ A=O$, $X\circ I=O$ and $[A,X]=O$.  
%Write $X=L+L\trans$, where $L$ is strictly lower triangular.  
%Then $O=[A,X]=[A,L]+[A,L\trans]$.  
%Since $[A,L]$ and $[A,L\trans]$ are strictly lower triangular and strictly upper triangular respectively, $[A,L]=[A,L\trans]=O$.  
%Also, since $M(P_n)=1$, all eigenvalues of $A$ are distinct and consequently $L$ is a polynomial in $A$.  
%This means $L$ is both symmetric and strictly lower triangular.  Hence $L=O$ and $X=O$, 
%and we conclude that $A$ has the SSP. By Theorem \ref{thm:SSP}, every supergraph of $P_n$ has a realization cospectral with $A$. 
%\end{ex}

\begin{obs}\label{exdiagnotSSP}
If the diagonal matrix $D= \diag(\mu_1, \mu_2, \ldots, \mu_n)$ has a repeated eigenvalue, 
say $\mu_i=\mu_j$, then $D$ does not have the SSP as validated by the matrix $X$ with a $1$ in positions $(i,j)$ and 
$(j,i)$ and zeros elsewhere.
\end{obs}

\begin{rem}\label{exdiagSSP}
Note that a diagonal matrix $D=\diag(\lambda_1,\lambda_2,\ldots, \lambda_n)$ with distinct eigenvalues has the SSP, because $DX=XD$ implies all off-diagonal entries of $X$ are zero.  Therefore, by Theorem \ref{thm:SSP}, every graph on $n$ vertices has a realization that is cospectral with $D$ and has the SSP.  
The existence of a cospectral matrix was proved in \cite{MS} via a different method.
\end{rem}

\noindent
However,  not every matrix with all eigenvalues distinct has the SSP.  

\begin{ex}\label{exdistinctnoSSP} Let 
\[A=\left(\begin{array}{rrrrr}
3 & -2 & 0 & 0 & 1 \\
-2 & 0 & -1 & 0 & 0 \\
0 & -1 & 3 & -1 & 0 \\
0 & 0 & -1 & 0 & 2 \\
1 & 0 & 0 & 2 & 3
\end{array}\right)\qquad\mbox{ and }\qquad X=
\left(\begin{array}{rrrrr}
0 & 0 & 1 & 1 & 0 \\
0 & 0 & 0 & 2 & -1 \\
1 & 0 & 0 & 0 & -1 \\
1 & 2 & 0 & 0 & 0 \\
0 & -1 & -1 & 0 & 0
\end{array}\right).  \]
Then $A$ does not have the SSP because 
  $A\circ X=O$, $I\circ X=O$, $[A,X]=O$, but $X\neq O$.  Note that $\spec(A)=\left\{2 \left(1+\sqrt{2}\right),4,\frac{1}{2} \left(1+\sqrt{17}\right),2 \left(1-\sqrt{2}\right),\frac{1}{2}
   \left(1-\sqrt{17}\right)\right\}.$
%Another example is the adjacency matrix of $K_{1,4}$.
\end{ex}

%%%%%%%%%%%%%%%%%%%%%%%%%%%%%%%%%%%%%%%%%%%%

\subsection{The Strong Multiplicity Property}
Let  $\oml=(m_1,\dots,m_q)$ be an ordered list of positive integers with $m_1 + m_2+ \cdots +  m_q=n$. 
We let $\mmu_{\oml}$  denote the set of all symmetric matrices whose ordered multiplicity list is $\oml$.
Thus, if $A$ has multiplicity list $\oml$, then 
\[
\mmu_{\oml}=\{B\in\Rsn: B \mbox{ has the same ordered multiplicity list as } A\}.
\]
It follows from results in \cite{Arnold} that $\mmu_{\oml}$ is a manifold. 
In the next lemma  we determine $\N_{\mmu_{\oml}. A}$.

\begin{lem}
\label{lem:SMP}
Let $A$ be an $n\times n$ symmetric matrix with exactly $q$ distinct eigenvalues, spectrum $\Lambda$, and ordered multiplicity list $\oml$.
Then 
\[
\N_{\mmu_{\oml}.A}= 
\{ X \in \mathcal{C}(A) \cap \Rsn : \tr(A^iX)=0 \mbox{ for } i=0,\ldots, q-1 \}.
\]
\end{lem}
\bpf
Let $A$ be an $n\times n$ symmetric matrix with spectrum given by the multiset $\Lambda$, and let $B(t)$ ($t\in (-1,1)$) 
be a differentiable path of matrices having the same ordered multiplicity list as $A$.
Let $A$  and $B(t)$ have spectral decomposition 
$A=\sum_{j=1}^q \lambda_j E_j$, 
and $B(t) = \sum_{j=1}^q \lambda_j(t) F_j(t)$, respectively. 
Then $\dot{B}(0)= \sum_{j=1}^q \dot{\lambda_j}(0) E_j + \sum_{j=1}^q \lambda_j \dot{F}_j(0)$, 
and we conclude that the tangent space  of    $\mmu_{\oml}$  is the sum of 
\[
\Tc:= \left\{ \sum_{j=1}^q c_jE_j: c_j \in \R \right\}
\]
  and  
the tangent space, $\T_{\mei_{\Lambda}.A}$. Therefore
\beq\label{neweq31c}
\N_{\mmu_{\oml}.A}= \N_{\mei_{\Lambda}.A} \cap \{S \in \Rsn : \tr(E_jS)=0 \mbox{ for all }j=1,\dots,q \}.
\eeq
As $\Span(E_1, \ldots, E_q )= \Span( I=A^0,  \ldots, A^{q-1})$, 
 $\tr(A^iS)=0$ for 
$i=0,  \ldots, q-1$ if and only if $\tr (E_jS)=0$ for  $j=1,\ldots, q$.   The result now follows. 
\epf

\begin{defn} The $n \times n$ symmetric matrix $A$ satisfies the {\it Strong Multiplicity Property} (or $A$ has the SMP for short)
provided the only symmetric matrix $X$ satisfying $A\circ X=O$, $I\circ X=O$, $[A,X]=O$ and $\tr(A^iX)=0$ for $i=0,\ldots, n-1$ is $X=O$. 
\end{defn}

\begin{rem}\label{cor:2evalSSPSMP} The minimal polynomial provides a dependency relation among the powers of $A$, so we can replace $i=0,\ldots, n-1$ by $i=0,\ldots, q-1$ where $q$ is the number of distinct eigenvalues of $A$.  
Since $I\circ X=O$ and $A\circ X=O$ imply $\tr(IX)=0$ and $\tr(AX)=0$, we can replace $i=0,\ldots, q-1$ by $i=2,\ldots, q-1$.  Therefore, for multiplicity lists with only $2$ distinct eigenvalues, the SSP and the SMP are equivalent.
\end{rem}

Lemma \ref{lem:SMP} asserts that $\mptn(G)$ and $\mmu_{\oml}$ intersect transversally at $A$ if and only if 
$A$ has the SMP.   A proof similar to that of
Theorem \ref{thm:SAP} yields the next result.
\begin{thm}
\label{thm:SMP}
If $A\in \mptn(G)$ has the SMP, then every supergraph of $G$ with the same vertex set has a realization that has the same ordered multiplicity list as $A$ and has the SMP.
\end{thm}

\begin{obs}
Clearly the  SSP implies the SMP. 
\end{obs}

 Example \ref{exdistinctnoSSP} and the next remark  show the SSP and the SMP are distinct.

\begin{rem}\label{distinctSMP}
In contrast to Example $\ref{exdistinctnoSSP}$,  every symmetric matrix whose eigenvalues are  distinct has the SMP.
To see this, let $A\in \Rsn$ have $n$ distinct eigenvalues and let its spectral decomposition be 
\[
A=\sum_{i=1}^n \lambda_j \by_j\by_j\trans.
\]
Suppose $X\in \Rsn$, $I\circ X=O$, $A\circ X=O$, $[A,X]=O$, and  $\tr(A^{i}X)=0$ for $i=0,1,\ldots, n-1$.
Since $A$ and $X$ commute, each eigenspace of $A$ is an invariant subspace of $X$.  These conditions  and the distinctness of eigenvalues imply that each $\by_j$ is an eigenvector of $X$, and    $\tr(p(A)X)=0$ for all polynomials $p(x)$. The distinctness of eigenvalues implies that $\by_j\by_j\trans$ is a polynomial in $A$ for $j=1,\dots,n$.
Hence $0=\tr(\by_j\by_j\trans X)= \by_j\trans X\by_j$, so   the eigenvalue of $X$ for which $\by_j$ is an eigenvector of $X$ is $0$ for $j=1,\dots,n$.
Thus  $X=O$,  and we conclude that $A$ has the SMP. 
\end{rem}

\begin{obs} Clearly the SMP implies the SAP.  
\end{obs}
The next example shows  that the SMP and the SAP are distinct.

\begin{ex}\label{SMPnotSAP} Consider the matrices
\[A=\left( \begin{array}{rrrrrrrr} 1 & 1 & 1 & 0 & 1 & 0 & 0 & 0 \\
 1 & 3 & 0 & 1 & 0 & 1 & 0 & 0 \\
 1 & 0 & 3 & -1 & 0 & 0 & 1 & 0 \\
 0 & 1 & -1 & 1 & 0 & 0 & 0 & 1 \\
 1 & 0 & 0 & 0 & 3 & -1 & -1 & 0 \\
 0 & 1 & 0 & 0 & -1 & 1 & 0 & -1 \\
 0 & 0 & 1 & 0 & -1 & 0 & 1 & 1 \\
 0 & 0 & 0 & 1 & 0 & -1 & 1 & 3
 \end{array} \right)\!, \
 X=\left(\begin{array}{rrrrrrrr} 
 0 & 0 & 0 & 1 & 0 & 0 & 1 & 0 \\
 0 & 0 & 1 & 0 & 0 & 0 & 0 & 1 \\
 0 & 1 & 0 & 0 & -1 & 0 & 0 & 0 \\
 1 & 0 & 0 & 0 & 0 & -1 & 0 & 0 \\
 0 & 0 & -1 & 0 & 0 & 0 & 0 & 1 \\
 0 & 0 & 0 & -1 & 0 & 0 & 1 & 0 \\
 1 & 0 & 0 & 0 & 0 & 1 & 0 & 0 \\
 0 & 1 & 0 & 0 & 1 & 0 & 0 & 0\end{array} \right)\!.\]
Clearly $A\circ X=O$ and $I\circ X=O$.  It is straightforward to verify (using computational software) that 
 $\spec(A)=\{0^{(4)}, 4^{(4)}\}$ and $[A,X]=O$, where $ \lambda^{(m)}$ denotes that eigenvalue $\lambda$ has multiplicity $m$. 
Since $A$ has only two eigenvalues, $A$ does not have the SMP by Remark \ref{cor:2evalSSPSMP}.  It is also straightforward to verify (using computational software) that both $A$ and $A-4I$ have the SAP.  Note that $\G(A)=Q_3$ and $A$ is a diagonal scaling of the positive semidefinite matrix  of rank  four constructed in \cite[Example 2.1]{Peters}. \end{ex}

\begin{obs}
If $A$ has the SSP (SMP), then $A + \lambda I$ has the SSP (SMP) for all $\lambda \in \mathbb{R}$.
\end{obs}

\begin{rem}
\label{exdistinctnoSSP2}
If $\lambda$ is the only multiple eigenvalue of $A$,
then $A$ has the {SMP} if and only if $A - \lambda I$ has the SAP. 
To see this assume that $\lambda_1$ is the only multiple eigenvalue of $A$, 
$A-\lambda_1 I$ has the SAP, 
$\lambda_2,  \ldots, \lambda_{q}$ are the remaining eigenvalues of $A$, and $\by_j$ is a unit eigenvector 
 of $A$ corresponding to $\lambda_j$ $(j=2, \ldots, q)$.
 Then the spectral decomposition of $A$ has the form 
  $A= \lambda_1 E_1 + \sum_{j=2}^q \lambda_j \by_j\by_j\trans$. 
  Assume $X$ is a symmetric matrix such that $A\circ X=O$, $I\circ X=O$, $[A,X]=O$, and  $\tr(A^kX)=0$
  for all $k$. As in Remark \ref{distinctSMP}, each $\by_j$ is an eigenvector of $X$ and $\by_j\by_j\trans$ is a polynomial in $A$, so $
  0=\tr(\by_j\by_j\trans X)= \by_j\trans X \by_j.$ Thus
   we conclude that $\by_j$ is in the null space of $X$ for $j=2,\ldots, q$. 
   Therefore, $AX=\lambda_1 E_1X=\lambda_1X$ (with the latter equality coming from $E_1+ \sum_{j=2}^q\by_j{\by}_j\trans=I$).
   Thus, $(A-\lambda_1 I)X=O$.  Since $A-\lambda_1 I$ has the SAP, $X=O$ and we conclude that $A$ has the SMP.
\end{rem}

%%%%%%%%%%%%%%%%%%%%%%%%%%%%%%%%%%%%%%%%%%%%
%%%%%%%%%%%%%%%%%%%%%%%%%%%%%%%%%%%%%%%%%%%%
\section{Properties of matrices having the SSP or SMP}\label{commonSSPSMP}

{Section \ref{sstrongman}  presents characterizations of the tangent spaces {$\T_{{\mrk_A}. A}$, $\T_{\mei_\Lambda.A}$, and $\T_{\mmu_{\bf m}.A}$} and applies these to obtain lower bounds on the number of edges in a graph where a matrix has the associated strong property.} Section \ref{sstrongcomp}  describes a computational test for determining whether a matrix has the SSP or the SMP.  Section \ref{ssGersh} presents  the Gershgorin intersection graph and uses it to test for the SSP.  Section \ref{sblockD} characterizes when block diagonal matrices have the SSP or the SMP in terms of the diagonal blocks.

%%%%%%%%%%%%%%%%%%%%%%%%%%%%%%%%%%%%%%%%%%%%

\subsection{Tangent spaces for the strong property manifolds}\label{sstrongman}
We begin by giving equivalent, but more useful, descriptions of the tangent spaces $\T_{{\mrk_A}. A}$, $\T_{\mei_\Lambda.A}$, and $\T_{\mmu_{\bf m}.A}$.  The set of all $n\times n$ skew-symmetric matrices is denoted by $\Rskewn$.

\begin{thm}
\label{thm:dim}
Let $A$ be an $n\times n$ symmetric matrix with {distinct eigenvalues $\lambda_1, \ldots, \lambda_q$}, spectral decomposition 
$A= \sum_{j=1}^q \lambda_j E_j$, ordered multiplicity list  ${\bf m}=(m_1, \ldots, m_q)$ and  rank $r${; the spectrum of $A$ is $\Lambda=\{\lambda_1^{(m_1)}, \ldots, \lambda_q^{(m_q)}\}$.}
Then 
\begin{itemize}
\item[\rm (a)] {\rm\cite{HLS}}
$\T_{{\mrk_A}. A} = \{AY+Y\trans A: Y  \mbox{ is an $n \times n$ matrix}\}$ and ${\dim}\; \T_{{\mrk_A}. A}= {n+1 \choose 2} -{n-r+1 \choose 2}= {r+1 \choose 2} + r(n-r)$;
\item[\rm (b)]
$\T_{\mei_\Lambda.A}= \{ AK-KA: K \in \Rskewn\}$ and ${\dim}\; \T_{\mei_\Lambda.A}= {n \choose 2} - \sum_{j=1}^q {m_i \choose 2}$; 
\item[\rm (c)]
$\T_{\mmu_{\bf m}.A}= \T_{\mei_{\Lambda .A}} + \Span \{ I=A^0, \ldots, A^{q-1}\}$ and ${\dim}\;  \T_{\mmu_{\bf m}.A} =  {n \choose 2} - \sum_{j=1}^q {m_i \choose 2}+ q$. \end{itemize}
\end{thm}

\bpf
Throughout the proof  $\{\bv_1, \ldots, \bv_n\}$ is an orthogonal basis of eigenvectors of $A$ with corresponding 
eigenvalues $\mu_1, \ldots, \mu_n$ with $\mu_i\neq 0$ for $i=1, \ldots, r$ and $\mu_i=0$ for $i=r+1, \ldots, n$.

To prove (a), consider an $n\times n$ matrix  $Y$ and $X \in \N_{\mrk_A .A}$.   First note that Lemma \ref{lem:SAP} asserts that   $X\in \Rsn$ and $AX=O$.  Thus  $XA=O$ and 
\[
\begin{array}{lll} 
\tr((AY+Y\trans A)X) & = & \tr( AYX+Y\trans AX)\\
& = & \tr( XAY+Y\trans AX)\\
& = & \tr(O+O)\\
& = & 0.
\end{array}
\]
{Therefore $\{AY+Y\trans A\}\subseteq \N_{\mrk_A .A}^{\perp}=\T_{{\mrk_A}. A}$.} Next note that  if $i\leq r$, then $\bv_i\bv_j\trans  +\bv_j\bv_i\trans = AY+Y\trans A$ where $Y= \frac{1}{\lambda_i} \bv_i\bv_j\trans $.
%{\gre [Deleted the `Similarly...' sentence because it is redundant]} % Similarly  if $j\leq r$, then $\bv_j\bv_i\trans + \bv_i\bv_j \trans= AY+Y\trans A$ where $Y=\frac{1}{\lambda_j} \bv_j\bv_i \trans$.
 Define $\Omega$ to be the set of pairs $(i,j)$  with {$1\leq i\leq j\leq n$ and  $i\leq r$.  Thus }
\begin{equation}
\label{eq:nice-tan-mrk}
\begin{array}{lll}
\Span \{ \bv_i\bv_j\trans  +\bv_j\bv_i\trans  : (i,j) \in \Omega\}  & \subseteq & \{ AY+ Y\trans A  : Y \mbox{ is an $n\times n$ matrix} \} \\
  &\subseteq & {\T_{\mrk_A. A}.}
 \end{array} 
 \end{equation}
   It is easy to verify that  $\{ \bv_i\bv_j\trans  +\bv_j\bv_i\trans  : (i,j) \in \Omega\}$ is an orthogonal set and hence
 \[ {\dim} \Span \{ \bv_i\bv_j\trans  +\bv_j\bv_i\trans  : (i,j) \in \Omega\}= {r+1 \choose 2}+ r(n-r). \]
 By \eqref{neweq31a}, $\dim \N_{\mrk_A .A}= {{n-r+1} \choose 2}$, and hence $\mbox{dim\;} \T_{\mrk_A.A}= {{n +1} \choose 2} - {{ n-r+1} \choose 2}$.
It follows that $ {\dim} \Span \{ \bv_i\bv_j\trans  +\bv_j\bv_i\trans  : (i,j) \in \Omega\}= \dim \T_{{\mrk_A}.A}$.  Therefore, equality holds throughout (\ref{eq:nice-tan-mrk}), 
and (a) holds.

To prove (b) consider  $K\in \Rskewn$ and $B \in \N_{\mei_{\Lambda}.A}$.   {Note} that  Lemma \ref{lem:SSP} asserts that  $B\in \mathcal{C}(A) \cap \Rsn$. Thus 
\[
\begin{array}{rcl}
\tr( (AK-KA)B)& = & \tr (AKB-KAB)\\
& = & \tr (K(BA-AB))\\
& = & \tr(O) \\
& = & 0.
\end{array}
\]
{Observe that if $\mu_i \neq \mu_j$, then $\bv_i\bv_j\trans +\bv_j\bv_i\trans = AK-KA$ where $K$ is the skew-symmetric matrix 
$\frac{1}{\mu_i-\mu_j}(\bv_i\bv_j\trans -\bv_j\bv_i\trans )$.}  Hence, 

\begin{equation}
\label{eq:nice-tan-mei}
\begin{array}{lll}
 \{ \bv_i\bv_j\trans  +\bv_j\bv_i\trans : \mu_i \neq \mu_j\}&  \subseteq&  \{ AK-KA: K \in \Rskewn \} \\
& \subseteq &  {\N_{\mei_{\Lambda} .A}}^{\perp}\\
&  =  & \T_{\mei_{\Lambda} .A}\\
& = & \{ \bv_i\bv_j\trans  +\bv_j\bv_i\trans : \mu_i \neq \mu_j\}.
\end{array}
\end{equation}
{where the last equality is \eqref{neweq31z}.  Thus, equality holds throughout (\ref{eq:nice-tan-mei}).  The dimension of $\Span \{ \bv_i\bv_j\trans  +\bv_j\bv_i\trans : {\mu_i \neq \mu_j}\}$
is easily seen to be ${n+1 \choose 2} - \sum_{j=1}^q {m_i+1 \choose 2}= {n \choose 2} -\sum_{j=1}^q  {m_i \choose 2}$,  and we have proven (b).}

Statement (c) follows from {(b)} and \eqref{neweq31c}, and the facts that \[\Span \{E_1,\ldots, E_q\} = \Span \{ I=A^0, \ldots, A^{q-1}\}\] and
\[\T_{\mei_{\Lambda}.A} \cap \Span \{ A^0, \ldots, A^{q-1}\} =\{O\}\]
since $\Span \{ A^0, \ldots, A^{q-1}\} \subseteq \mathcal{C}(A) \cap \Rsn= \N_{\mei_{\Lambda}.A}$.
\epf

\begin{rem} %{\gre [Changed from Observation to Remark and removed `$|E(G)|$ edges' and unwanted period in line 458]}
 \label{obs:nec}
Let $\Mc$ be a manifold in $\Rsn$ and $G$ be a graph of order $n$ %with $|E(G)|$ edges 
such that $\Mc$ and $\mptn(G)$
 intersect transversally at $A$.
 Then 
 \[{\dim}\;\T_{\Mc.A} + {\dim}\;\T_{\mptn(G).A} - {\dim}\; ( \T_{\Mc .A } \cap \T_{\mptn(G).A} ) \;{=\dim \Rsn}={ n+1 \choose 2}.\] 
Since  ${\dim}\;\T_{\mptn(G).A}= n + |E(G)|$  and ${ n+1 \choose 2}= {n\choose 2} +n$, 
\[
 |E(G)| = {n \choose 2} - {\dim}\; \T_{\Mc.A}+ {\dim} ( \T_{\Mc.A } \cap \T_{\mptn(G).A}).
\]
 \end{rem}
 
 \begin{cor}
 \label{cor:cond}
 {Let $G$ be a graph on $n$ vertices and let $A\in\mptn(G)$ with  spectrum $\Lambda$, ordered multiplicity list  ${\bf m}=(m_1, \ldots, m_q)$, and  rank $r$.} %, and let $|E(G)|$ be the number of edges in $G$.   
 Assume that $A$ is not a scalar matrix.
Then 
\begin{itemize}
\item[\rm (a)]{\rm\cite{HHMS10}} If $\mrk_{A}$ and $\mptn (G)$ intersect transversally at $A$, then 
\[ |E(G)| \geq  \left\{ \begin{array}{ll} { n-r+1 \choose 2} & \mbox{ if $G$ is not bipartite,}\\ [5pt]
  { n-r+1 \choose 2}-1 & \mbox{ if $G$ is bipartite;} \end{array} \right. 
\]
\item[\rm (b)]  If $\mei_{\Lambda}$ and $\mptn (G)$ intersect transversally at $A$, then $|E(G)| \geq \sum_{j=1}^q {m_i \choose 2}$; and
\item[\rm (c)]  If $\mmu_{\bf m}$ and $\mptn (G)$ intersect transversally at $A$, then $|E(G)| \geq \sum_{j=1}^q {m_i \choose 2} -  q+2$.
\end{itemize} 
 \end{cor}
 \bpf
Statement (b) follows immediately from {Remark \ref{obs:nec}} and part (b) of Theorem \ref{thm:dim}.

Statement (c) follows from {Remark \ref{obs:nec}}, part (c) of Theorem \ref{thm:dim}, and the fact that  $I$ and $A$ are linearly independent, and  lie in 
$\T_{\mmu_{\bf{m}.A}} \cap \T_{\mptn(G)}$. 

{Statement (a)   can be established using Remark \ref{obs:nec} and part (a) of Theorem \ref{thm:dim} but the argument is more complicated.  Since it was previously established in in \cite{HHMS10}  (see Theorem 6.5 and Corollary 6.6) we refer the reader there.}  \epf

%%%%%%%%%%%%%%%%%%%%%%%%%%%%%%%%%%%%%%%%%%%%

\subsection{Equivalent criteria for the strong properties}\label{sstrongcomp}

Let $H$ be a graph with vertex set $\{1,2,\ldots, n\}$ and edge-set $\{e_1, \ldots, e_p\}$, where $e_{k}= {i_kj_k}$.
For $A={(a_{i,j})} \in \Rsn$, we denote  the $p \times 1$ vector whose $k$th coordinate is $a_{i_k,j_k}$ by $\vect_H(A)$.
Thus, $\vect_H(A)$ makes a vector out of the elements of $A$ corresponding to the edges in $H$. 
Note that $\vect_H( \cdot )$ defines a linear transformation from $\Rsn$ to $\R^p$. {The {\em complement} $\overline{G}$ of $G$ is the graph with the same vertex set as G and edges exactly where G does not have edges.}

\begin{prop}
\label{prop:genman}
Let $\Mc$ be a manifold in $\Rsn$, let $G$ be a graph with vertices $1,2,\ldots, n$ such that $A \in \Mc \cap \mptn(G)$, and let 
$p$ be the number of edges in ${\overline{G}}$.
Then $\Mc$ and $\mptn(G)$ intersect transversally at $A$ if and only if $\{ \vect_{\overline{G}}\,(B): B\in \T_{\Mc .A} \} = \R^p$.
\end{prop}  
\bpf
First assume that $\Mc$ and $\mptn(G)$ intersect transversally at $A$. Consider ${\bc} \in \R^p$, and let ${C=(c_{ij})}\in\Rsn$  with ${c_{ij}=c_k}$ if 
${ij}$ is the $k$th edge of ${\overline{G}}$.   The assumption of transversality implies  there exist 
 ${B}\in \T_{\Mc .A}$ and ${D} \in \T_{\mptn(G).A}$ such that ${C=B+D}$.  Hence 
\[{\bc}=\vect_{\overline{G}}\,({C})= \vect_{\overline{G}}\,({B})+ \vect_{\overline{G}}\,(D)= \vect_{\overline{G}}\,(B),\] and 
we conclude that  $\{ \vect_{\overline{G}}\,(B): B\in \T_{\Mc .A} \} = \R^p$.

Conversely, assume that $\{ \vect_{\overline{G}}\,(B): B\in \T_{\Mc .A} \} = \R^p$ and consider $C\in \Rsn$. 
Then $\vect_{\overline{G}}\,(C)= \vect_{\overline{G}}\,(B)$ for some $B\in\T_{\Mc .A}$.  Note that  $C-B\in \T_{\mptn (G) . A}$.
Since  $C= B+ (C-B)$, we conclude that  $\T_{\Mc .A} + \T_{\mptn(G) . A}= \Rsn$.
\epf

In the {following, $E_{ij}$ denotes the $n\times n $ matrix with a $1$ in position $(i,j)$ and $0$ elsewhere, and
$K_{ij}$} denotes the $n\times n$ skew-symmetric matrix  $E_{ij}-E_{ji}$.
{Theorem \ref{thm:dim} %Corollary \ref{cor:cond}, 
and} Proposition \ref{prop:genman} imply the {next} result.

\begin{thm}
\label{thm:nec-suff-conds}{Let $G$ be a graph,  let $A\in\mptn(G)$ with
%Let $A\in \Rsn$ with 
$q$ distinct eigenvalues,} % and graph $G$,}
and let $p$ be the number of edges in ${\overline{G}}$. 
Then 
\begin{itemize}
\item[\rm (a)]   $A$ has the SAP if and only if the matrix whose columns are $\vect_{\overline{G}}\,(AE_{ij}+E_{ij}\trans A)$ for $1\leq i, j \leq n$ has rank $p$;
\item[\rm (b)]   $A$ has the SSP if and only if the matrix whose columns are   $\vect_{\overline{G}}\,(AK_{ij}-K_{ij} A)$ for $1\leq i < j \leq n$ has rank $p$; and
\item[\rm (c)] $A$ has the SMP if and only if the matrix whose columns are   $\vect_{\overline{G}}\,(AK_{ij}-K_{ij} A)$ for $1\leq i < j \leq n$ along with 
$\vect_{\overline{G}}\,(A^k)$ $(k=0,1,\ldots, q-1)$ has rank $p$.
\end{itemize}
\end{thm}

\begin{ex}
\label{bowtie} 
Let
\[
A=\left( \begin{array}{rrrrr}
1 & 1 & 0 & 0 & \sqrt{6} \\
1 & 1 & 0 &0 & \sqrt{6} \\
0 & 0 & 4 & 1 & 5\\
0 & 0 & 1 & 4 & 5 \\
\sqrt{6} & \sqrt{6} & 5 & 5 & 16
\end{array}
\right)\!\!.
\]
Then 
{\begin{eqnarray*}
[A,K_{1,3}]= \left(\begin{array}{rrrrr}
0 & 0 & -3 & -1 & -5 \\
 0 & 0 & 1 & 0 & 0 \\
 -3 & 1 & 0 & 0 & \sqrt{6} \\
 -1 & 0 & 0 & 0 & 0 \\
 -5 & 0 & \sqrt{6} & 0 & 0 \end{array} 
\right)\!\!, \quad  & 
[A,K_{1,4}]= \left( \begin{array}{rrrrr}
0 & 0 & -1 & -3 & -5 \\
 0 & 0 & 0 & 1 & 0 \\
 -1 & 0 & 0 & 0 & 0 \\
 -3 & 1 & 0 & 0 & \sqrt{6} \\
 -5 & 0 & 0 & \sqrt{6} & 0 \end{array} 
\right)\!\!, \\ [0pt]
[A,K_{2,3}]= \left( \begin{array}{rrrrr}
 0 & 0 & 1 & 0 & 0 \\
 0 & 0 & -3 & -1 & -5 \\
 1 & -3 & 0 & 0 & \sqrt{6} \\
 0 & -1 & 0 & 0 & 0 \\
 0 & -5 & \sqrt{6} & 0 & 0
 \end{array} 
\right)\!\!, \quad & 
[A,K_{2,4}] = \left( \begin{array}{rrrrr}
0 & 0 & 0 & 1 & 0 \\
 0 & 0 & -1 & -3 & -5 \\
 0 & -1 & 0 & 0 & 0 \\
 1 & -3 & 0 & 0 & \sqrt{6} \\
 0 & -5 & 0 & \sqrt{6} & 0
\end{array} \right)\!\!.
\end{eqnarray*}}
Let $G$ be the graph of $A$, and let $M$ be the matrix defined in part (b) of Theorem \ref{thm:nec-suff-conds} whose columns are  $\vect_{\overline{G}}\,([A,K_{ij}])$.  The submatrix $\widehat M$ of columns of $M$ corresponding to $K_{ij}$ where $(i,j)\in\{(1,3), (1,4), (2,3),(2,4)\}$   is
{\[\widehat M=
\left( \begin{array}{rrrr}
-3 & -1 & 1 & 0 \\
-1 & -3 & 0 &1\\
1 & 0 & -3 & -1 \\
0 & 1 & -1 & -3
\end{array} \right)\!\!.
\]}
Since $\widehat M$ is strictly diagonally dominant (equivalently, 0 is not in the union of Gershgorin discs of $\widehat M$),  $\widehat M$ is invertible, and so $\rank \widehat M=4$. %verified that this matrix has determinant 45, 
Therefore $M$ has rank 4 and by Theorem {\ref{thm:nec-suff-conds}}, we conclude that $A$ has the SSP.
\end{ex}

%{\gre [Former Cor 3.7 in Remark 2.19.]}
%\begin{cor}\label{cor:2evalSSPSMP} For multiplicity lists with only $2$ distinct eigenvalues, the SSP and the SMP are equivalent. \end{cor}

% \bpf Suppose that $A\in \Rsn$ has exactly two distinct eigenvalues, and let $G$ be the graph of $A$. The result follows from parts (b) and (c) of Corollary \ref{cor:cond}, and the facts that $q(A)=2$ and $\vect_{\overline{G}}\,(I)$ and $\vect_{\overline{G}}\,(A)$ are both the zero vector. \epf

%%%%%%%%%%%%%%%%%%%%%%%%%%%%%%%%%%%%%%%%%%%%

\subsection{Gershgorin discs and the SSP}\label{ssGersh}

Given a square matrix $A  = (a_{ij}) \in \mathbb{C}^{n \times n}$,
the \emph{Gershgorin intersection graph of $A$} is
the graph on vertices labeled $1, \dots, n$ in which
two vertices
$i \ne j$ are adjacent exactly when Gershgorin
discs $i$ and $j$ of $A$ intersect,
that is, when the inequality
\begin{equation}
\label{eqn:gersh}
\left|a_{ii} - a_{jj}\right| \mspace{8mu}
  \le \sum\limits_{\ell = 1, \ell \ne i}^n \mspace{-8mu}\left|a_{i\ell}\right|
  \mspace{12mu} +
  \sum\limits_{\ell = 1, \ell \ne j}^n \mspace{-8mu}\left|a_{j\ell}\right|
\end{equation}
is satisfied.
If $A$ has real spectrum, then
Gershgorin discs intersect if and only if they intersect on the real line,
and the Gershgorin intersection graph of $A$ is an interval graph.

Note that when graphs have a common vertex labeling, one of them
may be a subgraph up to isomorphism of another
without being identically a subgraph.  The next result requires the stronger condition of being identically a subgraph.  
\begin{thm}\label{thm:gersh}
Let $G$ be a graph with vertices labeled
$1, \dots, n$ and let  $A\in\mptn(G)$.
If the Gershgorin intersection graph of $A$
is identically a subgraph of $G$,
then $A$ satisfies the SSP.
\end{thm}
\bpf
Suppose that $e_k={i_{k}j_k}$ $(k=1,2, \ldots, p)$ are the edges of ${\overline{G}}$. 
Let $\widehat M$ be the $p \times p$ matrix whose $k$-th column 
is $ \vect_{\overline{G}}\,( AK_{i_{k} ,j_{k}}-K_{i_{k},j_{k}} A)$.

The $(k, k)$-entry of $\widehat{M}$ has  absolute value $\left| a_{i_{k}, i_{k}} - a_{j_{k}, j_{k} }\right|$
and the remaining  entries of  the $k$-th column of $\widehat{M}$ are, up to sign, a subset of the entries $a_{i_k,\ell}$, $\ell\ne i_k$, and $a_{j_k,\ell}$, $\ell \ne j_k$.
If the Gershgorin intersection graph of $A$ is identically a subgraph of $G$, 
then inequality (\ref{eqn:gersh}) is not satisfied for any $k$ (because the $e_k$ are nonedges of $G$ and therefore of the Gershgorin intersection graph).  Thus $\widehat M$ is strictly diagonally dominant, %so the union of the Gershgorin discs of $\widehat M$ does not contain $0$, 
so $\widehat M$ is invertible and has rank $p$.  Therefore, by Theorem \ref{thm:nec-suff-conds},  $A$ has the SSP.
\epf

Of course it is possible to have $\widehat M$ strictly diagonally dominant implying the invertibility of $\widehat M$ even when the Gershgorin intersection graph of $A$
is not  a subgraph of $G$, as in  Example \ref{bowtie}.

%%%%%%%%%%%%%%%%%%%%%%%%%%%%%%%%%%%%%%%%%%%%

\subsection{Block diagonal matrices}\label{sblockD}

\begin{thm}\label{dirsum} Let $A_i\in S_{n_i}(\R)$ for $i=1,2$. Then $A:=A_1\oplus A_2$ has the SSP (respectively, SMP)
 if and only if both $A_1$ and $A_2$ have the SSP (respectively, SMP) 
 and $\spec(A_1)\cap\spec(A_2)=\emptyset$.
\end{thm}
\bpf  Let $X=\mtx{X_1& W\\W\trans & X_2}$ be partitioned conformally with $A$. 

First, suppose  that $A_1$ and $A_2$ have the SSP, $\spec(A_1)\cap\,\spec(A_2)=\emptyset$,  $A\circ X=O, I\circ X = O$, and $[A,X]=O$.   Since $A_i\circ X_i=O$, $ I\circ X_i = O$, and $[A_i,X_i]=O$ for $i=1,2$ and $A_i$ has SSP, $X_i=O$ for $i=1,2$.  
The $1,2$-block of $[A,X]$ is $O=A_1W-WA_2$, so by \cite[Theorem 2.4.4.1]{HJ13}, $W=O$.  Thus $X=O$ and $A$ has the SSP.  

Now, suppose $A_1$ and $A_2$ have  the SMP rather than the SSP (and the spectra are disjoint).  Assume $A\circ X=O, I\circ X = O$,  $[A,X]=O$ and   $X$ also satisfies $\tr(A^kX)=0$ for $k=2,\dots,n-1$.  As before, $W=O$ so $X=X_1\oplus X_2$, $A_i\circ X_i=O$, $ I\circ X_i = O$,  $[A_i,X_i]=O$ for $i=1,2$. To obtain $X_i=O, i=1,2$ (and thus $X=O$ and $A$ has the SMP) it suffices to show that $\tr(A_i^kX_i)=0$ for $i=1,2$ and $k=2,\dots,n-1$.  
Consider the spectral decompositions of the diagonal blocks of $A$, 
$A_i=\sum_{j=1}^{q_i}\lambda_j^{(i)}E_j^{(i)}$
for $i=1,2$.  Since $\spec(A_1)\cap\spec(A_2)=\emptyset$, the spectral decomposition of $A$ is
\[A=\sum_{j=1}^{q_1}\lambda_j^{(1)}\left(E_j^{(1)}\oplus O\right)+\sum_{j=1}^{q_2}\lambda_j^{(2)}\left(O\oplus E_j^{(2)} \right)\!.\]
 Since each projection in the spectral decomposition is a polynomial in $A$, $A_1\oplus  O$ and $O\oplus A_2$ are polynomials in $A$.
Therefore,  $\tr(A_1^kX_1)=\tr\left(\left(A_1^k\oplus O\right)X\right)=0$ and $\tr(A_2^kX_2)=0$. %Since $A_1$ and $A_2$ have  the SMP,  $X_i=O, i=1,2$, so $X=O$ and $A$ has the SMP.
 
Conversely, assume $A_1\oplus A_2$ has the SSP and $A_1\circ X_1=O$, $ I\circ X_1 = O$, and $[A_1,X_1]=O$.  Then $X:=\mtx{X_1& O\\O & O}$ satisfies $A\circ X=O, I\circ X = O$, and $[A,X]=O$, so $X=O$, implying $X_1 =O$.  In the case of the SMP, $\tr (A_1^kX_1)=0$ implies $\tr (A^kX)=0$.  We show that $\spec(A_1)\cap\,\spec(A_2)\ne \emptyset$ implies $A$ does not have the SMP (and thus does not have the SSP):  Suppose $\lambda\in\spec(A_1)\cap\,\spec(A_2)$.  For $i=1,2$, choose $\bz_i\ne 0$ such that $A_i\bz_i=\lambda\bz_i$.  Define \[Z:=\mtx{O & \bz_1\bz_2\trans\\\bz_2\bz_1\trans & O}.\]
Then  $A\circ Z=O, I\circ Z = O$,  $\tr (A^kZ)=0$, and $[A,Z]=O$, but $Z\ne O$, showing that $A$ does not have the SMP. \epf

As one application we give an upper bound on $q(G)$ in terms of chromatic numbers.

\begin{thm}
Let $G$ be a graph and $\overline{G}$ its complement.  Then $q(G)\leq 2\chi(\overline{G})$.
\end{thm}
\bpf
The graph $G$ contains a disjoint union of $\chi(\overline{G})$ cliques.  Taking the direct sum of realizations of each clique each having at most two distinct eigenvalues and the SSP, and the eigenvalues of different cliques distinct gives a matrix having the SSP by Theorem \ref{dirsum}.  The result then follows from 
Theorem \ref{thm:SSP}.
\epf

Another application gives an upper bound on the number of distinct eigenvalues required for  a supergraph on a superset of vertices.

\begin{thm} 
\label{subgraph}
Let $A$ be a symmetric matrix of order $n$ with graph $G$.
If $A$ has the SSP (or the SMP) and $\widehat{G}$ is a graph on $m$ vertices containing $G$ as a subgraph, then there exists $\widehat{B}\in\mptn(\widehat{G})$ such that $\spec(A)\subseteq \spec(\widehat{B})$   (or has the  ordered multiplicity list of $A$ augmented with ones), and $\widehat B$ has the SSP (or the SMP).  Furthermore,    $q(\widehat{G}) \leq m -n+q(A)$.  If $A$ has the SSP, then we can prescribe $\spec(\widehat B)$ to be $\spec(A) \cup \Lambda$
where $\Lambda$ is any set of distinct real numbers such that $ \spec(A)\cap \Lambda=\emptyset$.
\end{thm}

\bpf
Assume that $A$ has the SSP (respectively, SMP), and without loss of generality 
that $V(G)=\{1,2,\ldots, n\}$,
$V(\widehat{G})=\{1,2,\ldots, m\}$, and $G$ is a subgraph of $\widehat{G}$.

Consider the matrix 
\[
B= \left(\begin{array}{cc} A& O \\
O & \diag(\Lambda) \end{array} \right)\!.
\]
Note that  the eigenvalues of $\diag(\Lambda)$ are distinct and distinct from the eigenvalues of $A$.    
It follows that  $\diag(\Lambda)$ has the SSP (see Remark
\ref{exdiagSSP} or note this follows from Theorem \ref{dirsum}), and thus has the SMP.    By Theorem \ref{dirsum}, if $A$ satisfies the SSP (SMP), then $B$ satisfies the SSP (SMP).

By Theorem \ref{thm:SSP} (or Theorem \ref{thm:SMP}), every supergraph of the graph of $B$ on the same vertex set has a realization $\widehat B$ that is cospectral with $B$ and has the SSP (or has the same ordered multiplicity list and has the SMP).
Hence $q(\widehat{G}) \leq q({B})=m-n+q(A)$.
\epf

\begin{rem}\label{rem:Fiedlr}  By taking a realization in $\mptn(G)$ with row sums $0$,
$\mr(G)\leq |G|-1$. It is a well known result that the eigenvalues 
of an irreducible tridiagonal matrix are distinct, that is $\mr(P_n)=n-1$.
 A classic result of Fiedler \cite{Fied} asserts that 
$\mr(G) =|G|-1$ if and only if $G$ is a path.

Theorem \ref{subgraph} can be used to derive this characterization, as follows.
  If $G$ contains a vertex of degree $3$ or more, 
then $G$ contains $K_{1,3}$ as a subgraph, and hence by Theorem \ref{subgraph} and Example \ref{exstar}
we conclude that $\mr(G) \leq |G|-2$. Also, it is easy to see if $G$ is disconnected,
then $\mr(G) \leq |G|-2$.  Thus, if $\mr(G) =|G|-1$, then $G$ has maximum degree $2$
and is connected.  Hence $G$ is a path or a cycle.  The adjacency matrix of a cycle $C$
has a  multiple eigenvalue, which implies that $\mr(C) \leq |C|-2$.
\end{rem}

%%%%%%%%%%%%%%%%%%%%%%%%%%%%%%%%%%%%%%%%%%%%
%%%%%%%%%%%%%%%%%%%%%%%%%%%%%%%%%%%%%%%%%%%%

\section{Application of strong properties to minimum number of distinct eigenvalues}\label{sq}

The SSP and the SMP allow us to characterize graphs having $q(G)\ge |G|-1$ (see Section \ref{shighq}).  First we introduce new parameters based on the minimum number of eigenvalues for matrices with the given graph that have the given strong property.

%%%%%%%%%%%%%%%%%%%%%%%%%%%%%%%%%%%%%%%%%%%%

\subsection{New parameters $\qs(G)$ and $\qm(G)$}\label{snewq}
%{\color{red}  Need to give definitions of $\xi$ and $\nu$.}
Recall that $\xi(G)$  is defined as the maximum nullity among matrices in $\mptn(G)$ that satisfy the SAP, and  $\nu(G)$ is the maximum nullity of positive semidefinite matrices having the SAP, so $\xi(G)\le\M(G)$ and $\nu(G)\le\Mplus(G)$.  These bounds are very useful because of the minor monotonicity of $\xi$ and $\nu$ (especially the  monotonicity on subgraphs). In analogy with these definitions, we define parameters for the minimum number of eigenvalues among matrices having the SSP or the SMP and described by a given graph.  In order to do this we need the property that every graph has at least one matrix with the property SSP (and hence SMP).  For any set  of $|G|$ distinct real numbers, there is matrix in $\mptn(G)$  with these eigenvalues that has the SSP by Remark \ref{exdiagSSP}.

\begin{defn}  We define
\[\qm(G)=\min\{q(A) :  A \in \mathcal{S}(G)\mbox{ and $A$ has the SMP}\}\]
and
\[\qs(G)=\min\{q(A) :  A \in \mathcal{S}(G)\mbox{ and $A$ has the SSP}\}.\]
\end{defn}

\begin{obs} From the definitions, $q(G)\le\qm(G)\le\qs(G)$ for any graph $G$.
\end{obs} 

One might ask why we have not defined a parameter $q_A(G)$ for the SAP.  The reason is that the SAP is not naturally associated with the minimum number of eigenvalues. The Strong Arnold Property considers only the eigenvalue zero; that is, if zero is a simple eigenvalue of $A$ (or not an eigenvalue of $A$), then $A$ automatically has the SAP. 
The next result is immediate from Remark \ref{cor:2evalSSPSMP} and the fact that $q(G)=1$ if and only if $G$ has no edges.
\begin{cor}
\label{cor:2evalqsqm}
Suppose $G$ is connected.  Then
$\qs(G)=2$ if and only if  $\qm(G)=2$.
\end{cor}

The next result is immediate from Theorem \ref{dirsum}.
\begin{cor}
\label{prop:dunion}
If $G$ is the disjoint union of connected components $G_i, i=1,\dots, h$ with $h\ge 2$, then
$\qs(G)=\sum_{i=1}^h \qs(G_i)$ and $\qm(G)=\sum_{i=1}^h \qm(G_i)$.
\end{cor}

\begin{rem}\label{rem:disjunion}  Suppose $G$ is the disjoint union of connected components $G_i, i=1,\dots, h$ with $h\ge 2$.  Since any graph has a realization %that has the SSP 
for any set of distinct eigenvalues  (Remark \ref{exdiagSSP}), $q(G)\le \max_{i=1}^h |G_i|$.  Clearly $q(G)\ge \max_{i=1}^h q(G_i)$.
%{For  the examples we have computed, $q(G)=\max_{i=1}^h q(G_i)$, but we do not know whether this equality holds for all $G$. } 
\end{rem}

The next result is an immediate corollary of Theorem \ref{subgraph}.
\begin{cor}
\label{cor:main}
If $H$ is a subgraph of $G$, $|H|=n$, and $|G|=m$, then
$q(G)\le \qs(G)\le m-n+\qs(H)$ and $q(G)\le \qm(G)\le m-n+\qm(H)$. \end{cor}

%%%%%%%%%%%%%%%%%%%%%%%%%%%%%%%%%%%%%%%%%%%%

\subsection{High values of $q(G)$}\label{shighq}

In this section we characterize graphs having $q(G)\ge |G| -1$.  

\begin{figure}[h]
\begin{center}
\scalebox{.7}{\begin{tikzpicture}
\draw[gray!0] (-2.5,-2.5) grid (2.5,1.5);
\node at (0,-3) {\Large 1) $H$ tree};
\node[whitenode] (v3) at (-1,0) {3};
\node[whitenode] (v4) at (1,0) {4};
\node (v1) [whitenode, above = 1 of v3] {1};
\node (v2) [whitenode, below = 1 of v3] {2};
\node (v5) [whitenode, above = 1 of v4] {5};
\node (v6) [whitenode, below = 1 of v4] {6};
\draw (v1)--(v3)--(v2);
\draw (v5)--(v4)--(v6);
\draw (v3)--(v4);
\end{tikzpicture}}
\scalebox{.7}{\begin{tikzpicture}
\draw[gray!0] (-2.5,-2.5) grid (2.5,1.5);
\node at (0,-3) {\Large 2) Campstool};
\node[whitenode] (v3) at (0,0) {3};
\foreach \x/\y in {135/1,45/2,225/4,-45/5}{
\node[whitenode] (v\y) at (\x:1.5) {\y};
\draw (v3)--(v\y);
}
%\draw (v1)--(v2);
%\end{tikzpicture}
%\begin{tikzpicture}
%\draw[gray!0] (-2.5,-2.5) grid (2.5,1.5);
%\node at (0,-3) {bowtie};
%\node[whitenode] (v5) at (0,0) {5};
%\foreach \x/\y in {150/1,210/2,30/3,-30/4}{
%\node[whitenode] (v\y) at (\x:1.5) {\y};
%\draw (v5)--(v\y);
%}
\draw (v1)--(v2);
\draw (v3)--(v4);
\end{tikzpicture}}
\\
\scalebox{.7}{\begin{tikzpicture}
\draw[gray!0] (-3,-1.5) grid (3,3);
\node at (0,-3.3) {\Large 3) Long $Y$ tree};
\node[whitenode] (v1) at (0,0) {1};
\foreach \x/\y/\z in {30/4/5,150/6/7,270/2/3}{
\node[whitenode] (v\y) at (\x:1.1) {\y};
\node[whitenode] (v\z) at (\x:2.2) {\z};
\draw (v\y)--(v\z);
\draw (v1)--(v\y);
}
\end{tikzpicture}}
\scalebox{.7}{\begin{tikzpicture}
\draw[gray!0] (-3,-2.5) grid (3,2);
\node at (0,-3.3) {\Large 4) $3$-sun};
\foreach \x/\y/\z in {270/1/4,30/2/5,150/3/6}{
\node[whitenode] (v\y) at (\x:1.1) {\y};
\node[whitenode] (v\z) at (\x:2.2) {\z};
\draw (v\y)--(v\z);
}
\draw (v1)--(v2)--(v3)--(v1);
\end{tikzpicture}}
\caption{Graphs for Proposition \ref{prop:HHY3}.}\vspace{-10pt}
\label{fig:HHY3}
\end{center}
\end{figure}

\begin{prop}
\label{prop:HHY3}
Let $G$ be one of  the  graphs  shown in Figure $\ref{fig:HHY3}$.  Then $\qs(G)\leq |G|-2$.
\end{prop}
\bpf
Let 
\[
A_1=\left(\begin{array}{rrrrrr}
0 & 0 & 1 & 0 & 0 & 0 \\
0 & 0 & 1 & 0 & 0 & 0 \\
1 & 1 & -1 & 1 & 0 & 0 \\
0 & 0 & 1 & 2 & 1 & 1 \\
0 & 0 & 0 & 1 & 1 & 0 \\
0 & 0 & 0 & 1 & 0 & 1
\end{array}\right)\!, \phantom{ and }
A_2=\left(\begin{array}{rrrrrr}
1 & 1 & 1 & 0 & 0 \\
 1 & 1 & -1 & 0 & 0 \\
 1 & -1 & 0 & 1 & 1 \\
 0 & 0 & 1 & 0 & 0 \\
 0 & 0 & 1 & 0 & 0
\end{array}\right)\!, 
%A_3=\left(\begin{array}{rrrrr}
%1 & 1 & 0 & 0 & \sqrt{6} \\
%1 & 1 & 0 & 0 & \sqrt{6} \\
%0 & 0 & 4 & 1 & 5 \\
%0 & 0 & 1 & 4 & 5 \\
%\sqrt{6} & \sqrt{6} & 5 & 5 & 16
%\end{array}\right)
\]
\[
A_3=\left(\begin{array}{rrrrrrr}
0 & 1 & 0 & 1 & 0 & 1 & 0 \\
1 & 1 & 1 & 0 & 0 & 0 & 0 \\
0 & 1 & 1 & 0 & 0 & 0 & 0 \\
1 & 0 & 0 & 1 & 1 & 0 & 0 \\
0 & 0 & 0 & 1 & 1 & 0 & 0 \\
1 & 0 & 0 & 0 & 0 & 1 & 1 \\
0 & 0 & 0 & 0 & 0 & 1 & 1
\end{array}\right)\!, \text{ and }
A_4= \left( \begin{array}{cccccc}
2 & 1 & 1 & 1 & 0 & 0\\
1 & 2 & 1 & 0 & 1 & 0 \\
1 & 1 & 2 & 0 & 0 & 1\\
1 & 0 & 0 & 1 & 0 & 0 \\
0 & 1 & 0 & 0 & 1 & 0\\
0 & 0 & 1 & 0 & 0 & 1
\end{array} \right)\!.
\]

The graphs of matrices $A_1,A_2,A_3$ and $ A_4$ are  the $H$ tree, the campstool, the long $Y$ tree, and the $3$-sun, respectively.  Also,  $q(A_i)= |G|-2$ for  $i=1,2,3,4$.
It is straightforward to  verify that each of the matrices $A_1,A_2,A_3$ and $ A_4$ has the SSP (see, for example, \cite{SageCode}).
\epf

\begin{cor}
\label{cor:nminus2}
If a graph $G$ contains a subgraph isomorphic to the  $H$ tree, the campstool, the long $Y$ tree,  or the $3$-sun
then \[q(G)\le q_S(G) \leq |G|-2.\]
\end{cor}

\bpf
This follows from Corollary \ref{cor:main} and Proposition \ref{prop:HHY3}.
\epf

%The maximum multiplicity $\M(G)$ is also the maximum nullity of a matrix in $\mptn(G)$.  The positive semidefinite maximum nullity $\Mplus(G)$ is the maximum nullity among positive semidefinite matrices in $\mptn(G)$) of a positive semidefinite  matrix in $\mptn(G)$.  
The parameters $\M(G)$ and $\Mplus(G)$ can be used to construct  lower bounds on $q(G)$.  %As $M(G)$ equals the maximum multiplicity among the eigenvalues over all matrices in $\mptn(G)$, 
For a graph $G$  of order $n$, clearly $q(G)\ge\lc\frac {n} {\M(G)}\rc$\!.  The next result improves on this in many cases, in particular for $G=K_{1,3}$.
\begin{prop}
\label{prop:Mplusbound}
For any graph $G$ on $n$ vertices,
% if $\Mplus(G)\le\lf\frac{n}{2}\rf$, then 
\[q(G)\geq 2+\lc\frac{n-2\Mplus(G)}{\M(G)}\rc\!.\]
Moreover,
if $\Mplus(G)<\frac {n} 2$, then $q(G)\ge 3$.
\end{prop}
\bpf
Let $A\in \mptn(G)$ be a matrix with $q(G)$ distinct eigenvalues.  Let $\alpha,\beta$ be the smallest and the largest eigenvalues of $A$.  Since $A-\alpha I$ and $-A+\beta I$ are positive semidefinite, the multiplicity of $\alpha$ and $\beta$ is no more than $\Mplus(G)$.  Every other eigenvalue of $A$ has multiplicity less than or equal to $\M(G)$.  Therefore,
% if $\Mplus(G)\le\lf\frac{n}{2}\rf$, the
$A$ has at least $2+\lc\frac{n-2\Mplus(G)}{\M(G)}\rc$ distinct eigenvalues.
The final statement of the proposition readily follows.
\epf

\begin{cor}
\label{cor:K3K13}
If $G$ contains two vertex disjoint subgraphs each of which is a $K_3$ or a $K_{1,3}$, 
then $q(G) \leq |G|-2$.
\end{cor}

\bpf
The $3\times 3$ all ones matrix $J_3$ has two distinct eigenvalues and has SSP, so $\qs(K_3)=2$.   Example \ref{exstar} 
and Proposition \ref{prop:Mplusbound} imply that  $\qs(K_{1,3})= 3$. Thus, by Corollary \ref{prop:dunion} and Corollary \ref{cor:main}, $q(G)\leq |G|-2$.
\epf

In Ferguson \cite[Theorem 4.3]{F80} it is shown that    a multiset of $n$ real numbers is the spectrum of an $n\x n$ periodic Jacobi matrix\footnote{A periodic Jacobi matrix is a real symmetric matrix with nonnegative off-diagonal entries whose graph is a cycle.} $A$  if and only if these numbers can be arranged as  $\lambda_1>\lambda_2\ge \lambda_3>\lambda_4\ge\lambda_5>\cdots\lambda_n$.  This solves the inverse eigenvalue problem for a cycle of odd length: Suppose $A\in\mptn(C_n)$  and $n$ is odd.  If the cycle product $a_{12}a_{23}\cdots a_{n-1,n}a_{n1}>0$, then $A$ is similar (by a diagonal matrix with diagonal entries in $\{\pm 1\}$) to a periodic Jacobi matrix;  if $a_{12}a_{23}\cdots a_{n-1,n}a_{n1}<0$, then $-A$ is similar  to a periodic Jacobi matrix.

 In the next result  we establish that a specific  matrix $A\in\mptn(C_n)$ with  $q(A)=\lc \frac n 2\rc$ has the SMP.
\begin{thm}
\label{thm:Ceimustar}
Let $C_n$ be the cycle on $n\ge 3$ vertices.  Then $\qm(C_n)=\lc\frac{n}{2}\rc$. \end{thm}

\bpf  Since $\M(C_n)=2$,  $q(C_n)\ge \lc \frac n 2\rc$.
Given $n \ge 3$, let $C=(c_{ij})$ be the {\em flipped-cycle matrix} of order $n$, that is, the (non-symmetric) $n\x n$ matrix
with entries $c_{i,i+1} = 1$ for $i=1,\dots,n-1$, $c_{n,1} = -1$, and $0$ otherwise.
Since $C$ satisfies the equation $C^n = -I$, the eigenvalues
of $C$ are the $n$th roots of $-1$.
The matrix $A = C+C\trans = C+C^{-1}$
is a symmetric matrix whose graph is the $n$-cycle $C_n$,
and whose eigenvalues
are $2\cos(2\pi\frac{2j-1}{2n})$ for $j \in \{ 1, \dots, n\}$,
which occur in $\lf\frac{n}2\rf$ pairs
satisfying $j_1+j_2 = n+1$, with
one singleton eigenvalue (coming from $2j=n+1$) equal to $-2$ when $n$ is odd.  Thus $q(A)=\lc \frac n 2\rc$.\vspace{3pt}

 We show that $A$ has the SMP, implying $q_M(C_n)=\lc \frac n 2\rc$.  Assume  $X=(x_{ij})$ is a symmetric matrix such that $A \circ X = O$, $I \circ X = O$,  $[A,X]=O$, % (or equivalently, $AX$ is symmetric), 
and $\tr(A^kX)=0$ for $k=1,\dots,n$. %; $[A,X]=O$ implies $A^kX$ is symmetric for $k=1,\dots,n$.  
Divide the  entries $x_{ij}$ (on or above the main diagonal)
 into $n$ bands for $k=0,\dots,n-1$ of the form $x_{i,i+k}$, $i=1,\dots,n-k$; all the entries of $X$ in bands $0,$ $ 1$ and $n-1$ are zero since  $I \circ X = O$ and $A \circ X = O$. %$\lceil \frac{n+1}2\rceil$  bands according to cyclic distance, i.e., the minimum distance from $i$ to $j$ mod $n$; in what follows, all arithmetic is done mod $n$.  
The fact that  $AX$ is symmetric implies that all the entries in each band are equal, and in addition,  $x_{1, 1+(n-k)}=-x_{1,1+k}$ for $k\le \frac n 2$. In the case that $n=2\ell$ is even, this implies $x_{i,i+\ell}=0$.   Now assume $X\ne O$ and let $m$ be the smallest natural number such that
band $m$ of $X$ contains a nonzero entry $x_{ij}$ (and $2\le m<\frac n 2$).
 Notice that the sign pattern $x_{i,i+m}=x_{1,1+m}$  for $i=1,\dots, n-m$ and $x_{i,i+(n-m)}=-x_{1,1+m}$ for $i= 1,\dots,m$ matches the sign pattern of $A^m$, i.e., $(A^m)_{i,i+m}=1$ for $i=1,\dots, n-m$ and  $(A^m)_{i,i+(n-m)}=-1$ for $i= 1,\dots,m$.  It is clear that $(A^m)_{ij}=0$ when the distance between $i$ and $j$ is greater than $m$, and by the choice of $m$, $x_{ij}=0$ when the distance between $i$ and $j$ is less than $m$.   Thus $\tr(A^mX)=2nx_{1,1+m}\ne 0$, a contradiction.  
  \epf

By the proof of Theorem \ref{thm:Ceimustar} and {Remark \ref{cor:2evalSSPSMP},}  the symmetric flipped cycle matrix $C+C\trans$ has the SSP for $n=4$.  It is clear from $X=C^2+{C^2}\trans$ that the symmetric flipped cycle matrix $C+C\trans$ does not have the SSP for $n\ge 5$ (this is implicit in  the  the proof of Theorem \ref{thm:Ceimustar}).
 
The next corollary is immediate from Corollary \ref{cor:main} and Theorem \ref{thm:Ceimustar}.
\begin{cor}
\label{cor:cycle}
If a graph $G$ contains a cycle of length $k$, then \[q(G)\le q_M(G)\leq |G|-\Big\lfloor\frac{k}{2}\Big\rfloor.\]
\end{cor}

\begin{prop}
\label{prop:qGequaln}
Let $G$ be a graph.  Then the following are equivalent:
\ben
\item[\rm(a)] $q(G)=|G|$,
\item[\rm(b)] $\M(G)=1$,
\item[\rm(c)] $G$ is a path.
\een
\end{prop}
\bpf
The equivalence of (b) and (c) is shown in \cite{Fied} (see also Remark \ref{rem:Fiedlr}). 
To show the equivalence of (a) and (b), first assume that $q(G)=|G|$. By Proposition 2.5 in \cite{AACFMN13}, $q(G)\leq \mr(G)+1$, so $\mr(G)=|G|-1$ and (b) holds.  Conversely, if $\M(G)=1$, then every eigenvalue has multiplicity at most one, so $q(G)=|G|$.  \epf

% add \label{rem:disjunion} to Remark 2.26
The next result characterizes all of the graphs $G$ that satisfy $q(G) \geq |G|-1$ and %, in fact, 
resolves a query presented in \cite{AACFMN13} on connected graphs that satisfy $q(G) = |G|-1$.

\begin{thm}\label{thm:qisGminus1}
{A graph $G$ has $q(G)\geq |G|-1$ if and only if}  $G$ is one of the following:
\ben
\item[\rm (a)] a path,
\item[\rm (b)]  the disjoint union of a path and an isolated vertex,
\item[\rm (c)]  {a path with one leaf attached to an interior vertex}, 
\item[\rm (d)] {a path  with an extra edge joining two vertices  at distance $2$}. 
\een
\end{thm}

\bpf
Suppose $G$ is a graph with $q(G)\geq |G|-1$.  {By  \cite[Proposition 2.5]{AACFMN13}, $q(G)\leq \mr(G)+1$, so  $\mr(G)\ge |G|-2$ and $\M(G)\le 2$.  }

 If $G$ has connected components {$G_i, i=1,\dots h$, with $h\ge 3$}  or at least two components containing two or more vertices, then {by Remark \ref{rem:disjunion}, $q(G)\leq \max |G_i|\le |G|-2$.}  If $G=H\dunion K_1$, then {$|H|=|G|-1=q(G)=q(H)$}, which  {implies $H$ is a path by Proposition \ref {prop:qGequaln}}.  So henceforth we assume $G$ is connected.

If $G$ has at least 2 cycles, then $G$ contains {as a subgraph at least one of}  $C_n$ with $n\ge 4$, two disjoint {copies of}  $K_3$, or the campstool.  By Corollaries \ref{cor:cycle}, \ref{cor:K3K13}, \ref{cor:nminus2}, it is impossible to have any of these graphs as a   subgraph.  Thus $G$ has one $C_3$ and no other cycles, or $G$ is a tree. 

{Suppose} that $G$ has one $3$-cycle.  Since the $3$-sun is forbidden {as a subgraph of $G$} by Corollary \ref{cor:nminus2}, at least one vertex {of the $C_3$ has degree two; let $v$ be this vertex.}  If $G -v$ is not a path, then $G-v$ has a vertex {$w$ of degree 3 or more.  Assume first that $w$ is not on the $C_3$ and note that  $w$ is a cut-vertex of $G$.  By choosing singular matrices for each component of $G-w$ and the matrix for the component containing $v$ having nullity two, we see that $\M(G-w)\geq 4$, leading to the contradiction that $\M(G)\ge 3$.    Now suppose $w$ is a vertex of $C_3$.  Then $G$ contains the campstool, contradicting Corollary \ref{cor:nminus2}.
Hence $G-v$ is a path, and $G$ consists of a path together with an extra edge joining two vertices at distance 2.}

Now suppose that $G$ is a tree.  Since  $\M(G)\leq 2$,  the maximum degree of $G$ is at most three by  a theorem of Parter and Wiener  (see, for example, \cite[Theorem 2.1]{FH}).  If the maximum degree is $2$, then $G$ is a path.  Otherwise, let {$w$} be a vertex {of  degree three}.  Then $\M(G)\leq 2$ implies that each {component of $G-w$} must be a path.  Since the $H$ tree is a forbidden subgraph by Corollary \ref{cor:nminus2}, the vertex of {each} path adjacent to {$w$} is an {endpoint}  of the path.  Since the long $Y$ tree is also forbidden by {Corollary} \ref{cor:nminus2}, at least one of the {components of $G-w$} is an isolated vertex, and we conclude that 
$G$ is path along with a pendent edge attached to an interior vertex of the path.

Each of the graphs $G$ listed  in the statement of the theorem has a unique path between two vertices of distance $|G|-2$, so by \cite[Theorem 3.2]{AACFMN13},  $q(G)\ge |G|-1$.  \epf

%%%%%%%%%%%%%%%%%%%%%%%%%%%%%%%%%%%%%%%%%%%%%

\section{Conclusions}
{We have presented a careful} development of the theory associated with the SSP and SMP as useful extensions of the SAP{.  A major consequence is that for any given  graph $G$ of order $n$ and multiset of $n$ real numbers (or partition of $n$), if one can find an SSP (or SMP) matrix $A\in\mptn(G)$ with spectrum equal to this multiset (or multiplicity list equal to this partition), then for each supergraph $\wt G$ on the same vertex set as $G$, that same spectrum (or multiplicity list) is attained by some matrix in $\mptn(\wt G) $.} Further implications to the inverse eigenvalue problem for graphs are expected, e.g., using these tools it may be possible to solve the inverse eigenvalue problem for additional families of graphs. A key issue moving forward, partly addressed here, is detection and construction of matrices in  $\mptn(G)$ that have the SSP or SMP. %We developed some useful tools along these lines, including  the use of Gershgorin's theorem and a discussion of block diagonal matrices. 

One of our main motivations  for considering the SSP and the SMP was for evaluating $q(G)$. In this work, we concentrated
on the case of graphs $G$ for which $q(G)$
was large relative to $|G|$. We established a number of forbidden subgraph-type results that were used to characterize the 
graphs $G$ such that $q(G) \geq |G|-1$. As a consequence, we produced a new verification of Fiedler's characterization that the path is the only graph $G$ such that $q(G)=|G|$, and we resolved an open problem left from the work in \cite{AACFMN13} concerning graphs $G$ with $q(G)=|G|-1$.
{A much clearer picture of $q(G)$ may be obtained by using the tools and techniques derived
here.}

%%%%%%%%%%%%%%%%%%%%%%%%%%%%%%%%%%%%%%%%%%%%%

\subsection*{Acknowledgments}

Shaun Fallat's research is supported in part by an NSERC Discovery Research Grant.  Leslie Hogben's research is supported in part by Holl Chair funds.

We also wish to thank Dr. Francesco Barioli for discussions with some of the authors that lead to the development and use of the Strong Spectral and Strong Multiplicity properties.

\providecommand{\noopsort}[1]{}


\begin{thebibliography}{1}


\bibitem{AACFMN13}
B.~Ahmadi, F.~Alinaghipour, M.S. Cavers, S.~M. Fallat, K.~Meagher, and
  S.~Nasserasr.
\newblock Minimum number of distinct eigenvalues of graphs.  
\newblock {\em Elec.~J.~Lin. Alg.}, 26:673--691, 2013.  % Erratum \url{http://repository.uwyo.edu/ela/vol26/iss1/45/}.

\bibitem{Arnold} V.I. Arnold.
On matrices depending on parameters.  In
{\em Vladimir I. Arnold - Collected Works},
Volume 2, pp 255-269,
 Springer-Verlag, Berlin, 2014.  %V.I. Arnold. On matrices depending on parameters. {\em Russian Math. Surveys}, 26:29--43, 1971.

\bibitem{Bar} F. Barioli. Private Communication, 2012.

\bibitem{param}	F. Barioli, W. Barrett, S. Fallat, H.T. Hall, 
L. Hogben, B. Shader, P. van den Driessche, and H. van der Holst.  
Parameters related to tree-width, zero forcing, and maximum 
nullity of a graph.  {\em J. Graph Theory}, 72: 146--177, 2013.

\bibitem{GCC}		F. Barioli, W. Barrett, S. Fallat, 
H.T. Hall, L. Hogben, and H. van der Holst. 
On the graph complement conjecture associated with minimum rank.
{\em Lin. Alg. Appl.}, 436: 4373--4391, 2012.

\bibitem{BF04}	F. Barioli and S.M. Fallat. 
 On two conjectures regarding an inverse eigenvalue problem for acyclic symmetric matrices.
{\em Elec.~J.~Lin. Alg.},
11: 41--50, 2004. 

\bibitem{BF05}	F. Barioli and S. M. Fallat.  On the eigenvalues of generalized stars and double generalized stars.
 {\em Lin. Multilin. Alg.},
53: 269--291, 2005.

\bibitem{BFH05}  F. Barioli, S.M. Fallat, and L. Hogben. A variant on the graph parameters of Colin de Verdi\`ere: Implications to the minimum rank of graphs. {\em Elec.~J.~Lin. Alg.}, 13: 387--404, 2005.    

\bibitem{BLMNSSSY13}	W. Barrett,  A. Lazenby, N. Malloy, C. Nelson, W. Sexton, R. Smith, J. Sinkovic, T. Yang.  The combinatorial inverse eigenvalue problem: complete graphs and small graphs with strict inequality.
{\em Elec.~J.~Lin. Alg.},
26: 656--672, 2013. 

\bibitem{CdV} 	Y. Colin de Verdi\`ere. 
\newblock On a new graph invariant and a criterion for planarity.
\newblock In {\em Graph Structure Theory}, pp. 137--147, 
American Mathematical Society, Providence, RI, 1993.

\bibitem{CdV2} 	Y. Colin de Verdi\`ere. 
\newblock Multiplicities of eigenvalues and tree-width graphs.
\newblock {\em J. Comb. Theory B},
74: 121--146, 1998.


\bibitem{DR}
A.~L.~Dontchev and R.~T.~Rockafeller. {\it Implicit Functions and Solution Mappings: A View form Variational Analysis}, 2nd ed., 
Springer Series in Operations Research and Financial Engineering, Springer,  2014.

\bibitem{FH}	S. Fallat and L. Hogben.
The minimum rank of symmetric matrices described by a graph: A survey.   
{\em Lin. Alg. Appl.}, 426: 558--582, 2007.  

\bibitem{HLA2}	S. Fallat and L. Hogben. Minimum Rank, Maximum Nullity, and Zero Forcing Number of Graphs.  In \emph{Handbook of Linear Algebra},  2nd ed., L. Hogben editor, CRC Press, Boca Raton, 2014. 

\bibitem{Fied}
M.~Fiedler.
\newblock A characterization of tridiagonal matrices.
 {\em Lin. Alg. Appl.}, 2:191--197, 1969.

\bibitem{F80}	W.E. Ferguson, Jr. The Construction of Jacobi and Periodic Jacobi
Matrices With Prescribed Spectra. {\em Math. Computation}, 35: 1203--1220, 1980.

\bibitem{HHMS10}
H.~T. Hall, L.~Hogben, R.~R. Martin, and B.~Shader.
\newblock Expected values of parameters associated with the minimum rank of a
  graph.
\newblock {\em Lin. Alg. Appl.}, 433:101--117, 2010.

\bibitem{HK} K.H.~Hoffman and R.~Kunze.
  {\em Linear Algebra},  2nd ed.
  Prentice Hall,  Upper Saddle River, NJ, 1971.

\bibitem{HLS}
H.~{\noopsort{Holst}}{van der Holst}, L.~Lov{\'a}sz, and A.~Schrijver.
\newblock The {C}olin de {V}erdi\`ere graph parameter.
\newblock In {\it Graph Theory and Computational Biology (Balatonlelle, 1996)},
  pp. 29--85, Janos Bolyai Math. Soc., Budapest, 1999.

\bibitem{HJ13}
R.~A. Horn and C.~R. Johnson.
\newblock {\em Matrix Analysis{\rm , 2nd ed}}.
\newblock Cambridge University Press, Cambridge, 2013.

\bibitem{JLS03} C.R. Johnson, A. Leal-Duarte, and C.M. Saiago.
    Inverse eigenvalue problems and lists of multiplicities of eigenvalues
    for matrices whose graph is a tree: the case of generalized stars and
    double generalized stars.
    \emph{Lin. Alg. Appl.} 373:311--330, 2003.

\bibitem{KP}
S.~G.~Krantz and H.~R.~Parks. {\em The Implicit Function Theorem: History, Theory and Applications}, Springer, 
2013.

\bibitem{L}
J.~M.~Lee. {\em Introduction to Smooth Manifolds}, 2nd ed., Graduate Texts in Mathematics, 
Springer, 2013.


\bibitem{MS}
K.~H.~Monfared and B.~L.~Shader. Construction of matrices with a given graph and prescribed interlaced spectral data.
{\em Lin. Alg. Appl.},  438: 4348-4358, 2015.

\bibitem{Peters}
T. Peters.
Positive semidefinite maximum nullity and zero forcing number.  
{\em Elec. J.  Lin. Alg.}, 23: 815--830,  2012.

\bibitem{SageCode} J.C.-H. Lin.  {\em Sage}  code for verifying the SSP.  Published on the Iowa State {\em Sage} server at \url{https://sage.math.iastate.edu/home/pub/52/}.
{\em Sage}  worksheet available at  \url{http://orion.math.iastate.edu/lhogben/has_SSP.sws}.

%
%%\bibliography{./JLaTeX/AuthorA,./JLaTeX/JournalA,./JLaTeX/JepBib}{}
%%\bibliographystyle{plain}
%
%\providecommand{\noopsort}[1]{}
%\begin{thebibliography}{1}
%
%\bibitem{AACFMN13Erratum}
%B.~Ahmadi, F.~Alinaghipour, M.S. Cavers, S.~M. Fallat, K.~Meagher, and
%  S.~Nasserasr.
%\newblock Erratum to ``number of distinct eigenvalues of graphs.
%\newblock {\em Electron.~J.~Linear Algebra}, 26:673--691, 2013.
%
%\bibitem{AACFMN13}
%B.~Ahmadi, F.~Alinaghipour, M.S. Cavers, S.~M. Fallat, K.~Meagher, and
%  S.~Nasserasr.
%\newblock Minimum number of distinct eigenvalues of graphs.
%\newblock {\em Electron.~J.~Linear Algebra}, 26:673--691, 2013.
%
%\bibitem{BFH05}  F. Barioli, S.M. Fallat, and L. Hogben. A variant on the graph parameters of Colin de Verdi\`ere: Implications to the minimum rank of graphs. {\em Elec. J.  Lin. Alg.}, 13: 387--404, 2005.    
%
%\bibitem{BF04}	F. Barioli and S.M. Fallat. 
% On two conjectures regarding an inverse eigenvalue problem for acyclic symmetric matrices.
% {\em Electron. J. Linear Algebra},
%11: 41--50, 2004. 
%
%\bibitem{FH}	S. Fallat and L. Hogben.
%The minimum rank of symmetric matrices described by a graph: A survey.   
%{\em Lin. Alg. Appl.} 426: 558--582, 2007.  
%
%\bibitem{HLA2}	S. Fallat and L. Hogben. Minimum Rank, Maximum Nullity, and Zero Forcing Number of Graphs.  In \emph{Handbook of Linear Algebra},  2nd edition, L. Hogben editor, CRC Press, Boca Raton, 2014. 
%
%\bibitem{Fied}
%M.~Fiedler.
%\newblock A characterization of tridiagonal matrices.
%\newblock {\em Linear Algebra Appl.}, 2:191--197, 1969.
%
%\bibitem{HHMS10}
%H.~T. Hall, L.~Hogben, R.~R. Martin, and B.~Shader.
%\newblock Expected values of parameters associated with the minimum rank of a
%  graph.
%\newblock {\em Linear Algebra Appl.}, 433:101--117, 2010.
%
%\bibitem{HLS}
%H.~{\noopsort{Holst}}{van der Holst}, L.~Lov{\'a}sz, and A.~Schrijver.
%\newblock The {C}olin de {V}erdi\`ere graph parameter.
%\newblock In {\it Graph Theory and Computational Biology (Balatonlelle, 1996)},
%  pp. 29--85, Janos Bolyai Math. Soc., Budapest, 1999.
%
%\bibitem{HJ13}
%R.~A. Horn and C.~R. Johnson.
%\newblock {\em Matrix Analysis{\rm , 2nd ed}}.
%\newblock Cambridge University Press, Cambridge, 2013.
%
%\bibitem{SageCode} J.C.-H. Lin.  {\em Sage}  code for verifying SSSP.  Published on the Iowa State {\em Sage} server at \url{https://sage.math.iastate.edu/home/pub/?/}.
%{\em Sage}  worksheet available at  \url{http://orion.math.iastate.edu/lhogben/?.sws}.
%
%\bibitem{Peters}
%T. Peters.
%Positive semidefinite maximum nullity and zero forcing number.  
%{\em Elec. J.  Lin. Alg.}, 23: 815--830,  2012.


\end{thebibliography}
\end{document}